\documentclass[11pt,twoside]{amsart}
\usepackage{amsfonts}
\usepackage{epsfig,graphics}
\usepackage{amssymb}
\usepackage{amscd}

\newtheorem{thm}{Theorem}[section]

\newtheorem{propn}[thm]{Proposition}
       
\newtheorem{defn}[thm]{Definition}

\newtheorem{lmm}[thm]{Lemma}
\newtheorem{rmk}[thm]{Remark}

\newtheorem{fact}[thm]{Fact}

\newcommand{\R}{{\bf R}}
\newcommand{\Z}{{\bf Z}}
\newcommand{\N}{{\bf N}}

\newcommand{\BP}{{\bf P}}
\newcommand{\X}{{\bf X}}

\newcommand{\aut}{\text{Aut }}

\newcommand{\rk}{\text{rk }}

\newcommand{\ad}{\text{ad }}

\newcommand{\SL}{\text{SL}}
\newcommand{\GL}{\text{GL}}
\newcommand{\PSL}{\text{PSL}}

\newcommand{\CO}{\text{CO}}

\newcommand{\Ein}{\text{Ein}}
\newcommand{\Conf}{\text{Conf }}

\newcommand{\Ad}{\text{Ad }}
\newcommand{\Aut}{\text{Aut }}

\newcommand{\liea}{{\mathfrak{a}}}

\newcommand{\lien}{{\mathfrak{n}}}
\newcommand{\lieu}{{\mathfrak{u}}}
\newcommand{\g}{{\mathfrak{g}}}
\newcommand{\h}{{\mathfrak{h}}}
\newcommand{\p}{{\mathfrak{p}}}
\newcommand{\s}{{\mathfrak{s}}}
\newcommand{\q}{{\mathfrak{q}}}
\newcommand{\lier}{{\mathfrak{r}}}
\newcommand{\oo}{{\mathfrak{o}}}
\newcommand{\co}{{\mathfrak{c}\mathfrak{o}}}
\newcommand{\gl}{{\mathfrak{g}\mathfrak{l}}}

\newcommand{\tx}{\tilde{x }}

\newcommand{\tz}{\tilde{z}}
\newcommand{\tp}{\tilde{p}}

\newcommand{\tdelta}{\tilde{\Delta }}
\newcommand{\tlambda}{\tilde{\Lambda}}

\newcommand{\tein}{\widetilde{\Ein }^{p,q}}
\newcommand{\tb}{\tilde{b }}

\newenvironment{Pf}{\medskip \noindent {\bf Proof: }}
   {$\diamondsuit$ }

\begin{document}
\pagenumbering{arabic}

\title[Conformal actions of nilpotent groups]{Conformal actions of nilpotent groups on pseudo-Riemannian manifolds}
\author{Charles Frances and Karin Melnick}
\date{\today}
\maketitle

\begin{abstract}
We study conformal actions of connected nilpotent Lie groups on compact
pseudo-Riemannian manifolds. We prove that if a type-$(p,q)$ compact manifold $M$ supports a conformal action of a connected nilpotent group $H$, then the degree
of nilpotence of $H$ is at most $2p+1$, assuming $p \leq q$; further, if this maximal
degree is attained, then $M$ is conformally equivalent to the universal
type-$(p,q)$, compact, conformally flat space, up to finite or cyclic covers. The proofs
make use of the canonical Cartan geometry associated to a pseudo-Riemannian
conformal structure.
\end{abstract}

\section{Introduction}

Let $(M,\sigma)$ be a compact \emph{pseudo-Riemannian} manifold---that is, the
tangent bundle of $M$ is endowed with a type-$(p,q)$ inner product, where
$p+q=n= \dim M$. We will always assume $p \leq q$.  The \emph{conformal class} of $\sigma$ is 
$$[ \sigma ] = \{ e^h \sigma \ : \ h : M \rightarrow \R \ \mbox{smooth} \}$$  Denote by $\Conf M$ the group of \emph{conformal
  automorphisms} of $M$---the group of diffeomorphisms $f$ of $M$ such that
$f^* \sigma \in [ \sigma ]$.
If $n \geq 3$, then $\Conf M$ endowed with the compact-open topology  is a Lie group (see \cite[IV.6.1]{kob} for the Riemannian case; the proof is similar for $p > 0$).

A basic question, first addressed by A. Lichnerowicz, is to characterize the pseudo-Riemannian manifolds $(M,\sigma)$ for which $\Conf M$ does not preserve any metric in $[ \sigma ]$; in this case, $\Conf M$ is {\it essential}.  The \emph{pseudo-Riemannian Lichnerowicz conjecture} says that if $M$ is compact and $\Conf M$ is essential, then $M$ is \emph{conformally flat}---that is, locally conformally equivalent to $\R^{p,q}$ with the translation-invariant, type-$(p,q)$ metric.  A stronger result in the Riemannian case was proved by Lelong-Ferrand \cite{lf1}.  See for example \cite{groupe-essentiel} for more background on this conjecture.




One difficulty for general type $(p,q)$ is that no characterization of essential conformal groups exists.  In the Riemannian case, on the other hand, for $M$ compact, $\Conf M$ is essential if and only if it is noncompact.  For $p \geq 1$, noncompactness is only a necessary condition to be essential.  Now, a first approach to the conjecture is to exhibit  sufficient conditions on a group of conformal transformations which ensure it is essential, and to test the conjecture on groups satisfying the given condition.  

From the conformal point of view, the \emph{Einstein universe} is the most symmetric structure of type $(p,q)$.  These spaces, denoted $\Ein^{p,q}$, are defined in section \ref{subsection.ein.structure} below.  The group $\Conf \Ein^{p,q}$ is isomorphic to $\mbox{PO}(p+1,q+1)$, and it is essential.  
The Einstein spaces are conformally flat.

Thanks to \cite{zirank}, we know that a simple noncompact real Lie group acting isometrically on a compact pseudo-Riemannian manifold $(M,\sigma)$ of type $(p,q)$ satisfies $\rk H \leq p$, where $\rk H$ denotes the real rank.  For $H < \Conf M$ noncompact and  simple, the rank
$$ \rk H \leq p+1=\rk \mbox{PO}(p+1,q+1)$$

This was first proved in \cite{zirank}, also in \cite{bn}, and for $H$ not
necessarily simple in {\cite[1.3 (1)]{bfm}}.  Thus, conformal actions of simple groups $H$, with $\rk H=p+1$, on type-$(p,q)$ compact pseudo-Riemannian manifolds cannot preserve any metric in the conformal class. The results of \cite{bn}, together with
\cite{fz}, give that when $H < \Conf M$ attains this maximal rank, then $M$ is
globally conformally equivalent to $\Ein^{p,q}$, up to finite covers when $p
\geq 2$; for $p=1$, $M$ is conformally equivalent to
the universal cover $\widetilde{\Ein}^{1,n-1}$, up to cyclic and finite
covers. In particular, $M$ is conformally flat, so this result supports the pseudo-Riemannian Lichnerowicz conjecture.
The interested reader can find a wide generalization of this theorem in {\cite[1.5]{bfm}}.

Actions of semisimple Lie groups often exhibit rigid behavior partly because the algebraic structure of such groups is itself rigid.  The structure of nilpotent Lie groups, on the other hand, is not that well understood; in fact, a classification of nilpotent Lie algebras is available only for small dimensions.  From this point of view, it seems challenging to obtain global results similar to those above for actions of nilpotent Lie groups.  Observe also that a pseudo-Riemannian conformal structure does not naturally define a volume form, so that the nice tools coming from ergodic theory are not available here.

 For a Lie algebra $\h$, we adopt  the notation $\h_{1} = [\h,\h]$, and $\h_{k}$ is defined inductively as $[\h,\h_{k-1}]$.  The \emph{degree of nilpotence} $d(\h)$ is the minimal $k$ such that $\h_{k} = 0$. For a connected, nilpotent Lie group $H$, define the nilpotence degree $d(H)$ to be $d(\h)$.  If a connected Lie group $H$ is  nilpotent 
and acts isometrically on a type-$(p,q)$ compact pseudo-Riemannian manifold $M$, where
$p \geq 1$, then
the nilpotence degree $d(H) \leq 2p$ (when $p=0$, then $d(H) =1$).  This was proved in the Lorentzian case
in \cite{zilor}, and in broad generality in {\cite[1.3 (2)]{bfm}}.  Theorem 1.3 (2) of  \cite{bfm} also implies $d(H) \leq 2p+2$ for $H < \Conf M$.  This bound is actually not tight, and the first result of the paper is to provide the tight bound, which turns out to be $2p+1$, the maximal nilpotence degree of a connected nilpotent subgroup in $\mbox{PO}(p+1,q+1)$.

\begin{thm}
\label{degree.bound}
Let $H$ be a connected nilpotent Lie group acting conformally on a  compact
pseudo-Riemannian manifold $M$ of type $(p,q)$, where $p \geq 1$, $p+q \geq 3$.  Then $d(H) \leq 2p+1$.
\end{thm}

By theorem 1.3 (2) of \cite{bfm}, a connected nilpotent group $H$ such that $d(H)=2p+1$ cannot act isometrically on a compact pseudo-Riemannian manifold of type $(p,q)$.  The following theorem says that if this maximal nilpotence degree is attained in
$\Conf M$, then $M$ is a complete conformally flat manifold, providing further support
for the pseudo-Riemannian Lichnerowicz conjecture.

\begin{thm}
\label{flatness.thm}
Let $H$ be a connected nilpotent Lie group acting conformally on a compact
pseudo-Riemannian manifold $M$ of type $(p,q)$, with $p \geq 1$, $p+q \geq 3$.  If $d(H) = 2p+1$, then $M$ is conformally equivalent to ${\widetilde \Ein}^{p,q}/ \Phi$, where $\Phi  < \widetilde{\mbox O}(p+1,q+1)$ is trivial or isomorphic to $\Z_2$ when $p \geq 2$, or isomorphic to $\Z$ when $p=1$. 
\end{thm}

Here, ${\widetilde \Ein}^{p,q}$ denotes the universal cover of $\Ein^{p,q}$ and $\widetilde{\mbox{O}}(p+1,q+1) = \Conf \tein$.  When $p \geq 2$, the center of $\widetilde{\mbox O}(p+1,q+1)$ has order two, while for $\widetilde{\mbox O}(2,q+1)$, $q \geq 2$, the center is infinite cyclic.  From this theorem follows a complete description of the $H$-action: let $\Phi'$ be the intersection of $\Phi$ with the center of $\widetilde{\mbox{O}}(p+1,q+1)$.  Then the conformal group of the quotient $\tein / \Phi$ is $N / \Phi'$, where $N$ is the normalizer of $\Phi$ in $\widetilde{\mbox{O}}(p+1,q+1)$.  The conformal diffeomorphism given by the theorem conjugates the $H$-action on $M$ to one on $\tein / \Phi$, via a representation $H \rightarrow N / \Phi'$, which is faithful if the original $H$-action on $M$ is.

We do not treat the Riemannian case $p=0$ in this paper, since the results in this case are a trivial consequence of Ferrand's theorem.

Section \ref{section.geometrique} below provides brief background on the geometry of the Einstein universe, as well as an algebraic study of nilpotent subalgebras of $\oo(p+1,q+1)$.  In section \ref{section.cartan}, we introduce the notion of Cartan geometry, which is central in all the proofs, and recall the interpretation of type-$(p,q)$ conformal structures, where $p+q \geq 3$, as Cartan geometries infinitesimally modeled on $\Ein^{p,q}$.  Section \ref{section.general.degree.bound} uses results of \cite{bfm} to prove theorem \ref{degree.bound}.  Actually, we prove here a stronger statement, theorem \ref{degree.bound.precise}, which gives also the starting point to prove theorem \ref{flatness.thm}: whenever a connected nilpotent group of maximal nilpotence degree acts conformally on $M$, then some point has nontrivial stabilizers, containing special elements called \emph{lightlike translations}.  We explain the role of these elements and outline the proof of theorem \ref{flatness.thm} at the end of section \ref{section.general.degree.bound}.

\section{${\Ein}^{p,q}$ as a homogeneous space for $\mbox{PO}(p+1,q+1)$}
\label{section.geometrique}
In this section, we introduce some basic notation used throughout the paper, and provide background on the geometry of the Einstein universe, as well as an algebraic study of nilpotent subalgebras of ${\mathfrak o}(p+1,q+1)$.
\subsection{Geometry of $\Ein^{p,q}$}
\label{subsection.ein.structure}
Let $\R^{p+1,q+1}$ be the space $\R^{p+q+2}$ endowed with the quadratic form
$$ Q^{p+1,q+1}(x_0, \ldots, x_{n+1}) = 2(x_0x_{p+q+1} + \cdots + x_{p}x_{q+1}) + \Sigma_{p+1}^q x_i^2$$
 We consider the null cone  
$$\mathcal{N}^{p+1,q+1} = \{  x \in \R^{p+1,q+1} \ | \  Q^{p+1,q+1}(x)=0 \}$$  
and denote by $\mathcal{\widehat{N}}^{p+1,q+1}$ the cone  $\mathcal{N}^{p+1,q+1}$ with the origin removed.
The projectivization ${\BP}(\widehat{\mathcal{N}}^{p+1,q+1})$ is a smooth submanifold of $\R\BP^{p+q+1}$, and inherits from the pseudo-Riemannian structure of  $\R^{p+1,q+1}$ a type-$(p,q)$ conformal class (more details can be found in \cite{charlesthese}, \cite{primer}).  We call the {\it Einstein universe} of type $(p,q)$, denoted $\Ein^{p,q}$, this compact manifold  ${\BP}(\widehat{\mathcal{N}}^{p+1,q+1})$ with this conformal structure. 
Note that $\Ein^{0,q}$ is conformally equivalent to the round sphere $({\bf S}^q,g_{{\bf S}^q}).$  When $p \geq 1$, the product $({\bf S}^p \times {\bf S}^{q},-g_{{\bf S}^p} \oplus g_{{\bf S}^q})$ is a conformal double cover of $\Ein^{p,q}$. 

The projective orthogonal group of $Q^{p+1,q+1}$, isomorphic to $\mbox{PO}(p+1,q+1)$,  acts projectively on $\Ein^{p,q}$ and is the full conformal group of $\Ein^{p,q}$.

\subsubsection{Lightcones, stereographic projection, and Minkowski charts}
\label{intuition}
A {\it lightlike}, {\it timelike}, or {\it spacelike} curve of a pseudo-Riemannian manifold $(M, \sigma)$ is a $C^1$ $\gamma : I \to M$ such that $\sigma_{\gamma(t)}(\gamma^{\prime}(t),\gamma^{\prime}(t))$ is $0$, negative, or positive, respectively, for all $t \in I$.  It is clear that the notion of lightlike, timelike and spacelike curves is a conformal one.  Lightlike curves are sometimes also called \emph{null}.

It is  a remarkable fact  that all metrics in $[\sigma]$ have the same null geodesics, as unparametrized curves (see  for example \cite{liouville} for a proof). Thus it makes sense to speak of null---or lightlike---geodesics for pseudo-Riemannian conformal structures. Given a point $x \in M$, the {\it lightcone} of $x$, denoted $C(x)$, is the set of all lightlike geodesics passing through $x$. 

\begin{figure}[ht]
 \includegraphics[width=2.5in,height=1.25in]{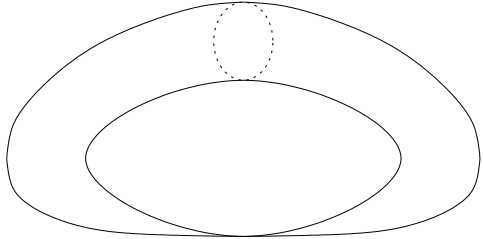}
\caption{the lightcone of a point in $\mbox{Ein}^{1,2}$}
\protect{\label{figure.lightcone}}
\end{figure}

The lightlike geodesics of $\Ein^{p,q}$ are the projections on $\Ein^{p,q}$ of totally isotropic $2$-planes in $\R^{p+1,q+1}$.  Hence every null geodesic is closed.
 If $x \in \Ein^{p,q}$ is the projection of $y \in {\mathcal N}^{p+1,q+1}$, the lightcone $C(x)$ is just $\BP({ y}^{\bot} \cap \mathcal{N}^{p+1,q+1})$.  Such a lightcone is not smooth, but  $C(x) \setminus \{x \}$ is smooth and diffeomorphic to $\R \times {\bf S}^{p-1}\times{\bf S}^{q-1}$ (see figure ~\ref{figure.lightcone}).  

For any $(p,q)$, there is a  generalized notion of  stereographic projection.  Consider $\varphi : \R^{p,q} \to \Ein^{p,q}$ given in projective coordinates of $\R \BP^{n+2}$ by
$$\varphi: x \mapsto [-\frac{1}{2}Q^{p,q}(x,x): x_1 : \cdots : x_n : 1]$$ 
Then $\varphi$ is a conformal embedding of $\R^{p,q}$ into $\Ein^{p,q}$, called the inverse {\it stereographic projection} with respect to $[e_0]$.  
The image $\varphi(\R^{p,q})$ is a dense open set of $\Ein^{p,q}$ with boundary the lightcone $C([e_0])$.  Since the action of $\mbox{PO}(p+1,q+1)$ is transitive on $\Ein^{p,q}$, it is clear that the complement of any lightcone $C(x)$ in $\Ein^{p,q}$ is conformally equivalent to $\R^{p,q}$. Such an open subset of $\Ein^{p,q}$ will be called a {\it Minkowski component}, and denoted ${\bf M}(x)$. Its identification with $\R^{p,q}$ via stereographic projection with respect to $x$ is a {\it Minkowski chart}.


Here we explain how images of lightlike lines of $\R^{p,q}$ under $\varphi$ reach the boundary (see also \cite[ch 4]{charlesthese}).  Lightlike lines of $\R^{p,q}$ are identified via $\varphi$ with traces on ${\bf M}([e_0])$ of lightlike geodesics in $\Ein^{p,q}$.  If $\gamma : \R \to \R^{p,q}$ is a lightlike line, then $\lim_{t \to \infty}\varphi(\gamma(t))=\lim_{t \to -\infty}\varphi(\gamma(t))=x_{\gamma}$, where $x_{\gamma} \in C([e_0]) \backslash \{ [e_0] \}$.  
For lightlike lines $\gamma(t) = c + tu$ and $\beta(t) = b + tv$, the limits $x_{\gamma}=x_{\beta}$ if and only if 
$u = v$ and $\langle b-c,u \rangle = 0$.  In other words, ${\bf M}([e_0]) \cap C(x)$ is a degenerate affine hyperplane for $x \in C([e_0]) \setminus \{ [e_0] \}$.

\subsubsection{A brief description of $\oo(p+1,q+1)$}  
\label{brief.description}
The Lie algebra $\oo(p+1,q+1)$ consists of all $(n+2) \times (n+2)$ matrices $X$, where $n=p+q$, such that 
$$ X^t J_{p+1,q+1} + J_{p+1,q+1} X = 0$$
where $J_{p+1,q+1}$ is the matrix of the quadratic form $Q^{p+1,q+1}$.
It can be written as a sum $\lieu^- \oplus \lier \oplus \lieu^+$ (see \cite[IV.4.2]{kob} for $p=0$; the case $p>0$ is a straightforward generalization), where

$$ \lier =  \left\{ \left( \begin{array}{ccc}
a &  & 0 \\
  & M   &  \\
  &    & -a
\end{array} \right) \ :    
\qquad 
\begin{array}{c}
 a \in \R  \\
  M \in \oo(p,q) \\
    \end{array} 
\right\}
$$

$$ \lieu^+= \left\{ \left( \begin{array}{ccc}
0& -x^t.J_{p,q}   &  0\\
  & 0  & x \\
  &   & 0
\end{array} \right) \ :  
\qquad 
\begin{array}{c}
  x\in \R^{p,q} 
\end{array}
\right\}
$$


and

$$ \lieu^- = \left\{ \left( \begin{array}{ccc}
0&    &  \\
 x & 0  &  \\
 0 & -x^t.J_{p,q}  & 0
\end{array} \right) \ :  
\qquad 
\begin{array}{c}
  x\in \R^{p,q} 
\end{array}
\right\}
$$

 Thus $\lier \cong \co(p,q)$, and there are two obvious isomorphisms $i^+$ and $i^-$ from $\lieu^+$ and  $\lieu^-$, respectively, to $\R^{p,q}$, given by the matrix expressions above.  


The standard basis of $\R^{p,q}$ corresponds under $i^-$ to the basis of $\lieu^-$

$$
U_i  =  \left\{ \begin{array}{cl} 
E_i^0 - E_{n+1}^{n+1-i} & i \in \{ 1, \ldots, p \} \cup \{ q+1, \ldots, n \}  \\
 E_i^0 - E_{n+1}^i & i \in \{ p+1, \ldots, q \}
\end{array} \right.
$$

where $E_i^j$ is the $(n+2)$-dimensional square matrix with all entries $0$ except for a $1$ in the $(i,j)$ place.

The parabolic Lie algebra $\p \cong \lier \ltimes \lieu^+$ is the Lie algebra of the stabilizer $P$ of $[e_0]$ in $\mbox{PO}(p+1,q+1)$, and similarly for $\p^- \cong \lier \ltimes \lieu^-$, the Lie algebra of the stabilizer of $[e_{n+1}]$.  The groups $P$ and  $P^-$ are isomorphic to the semidirect product $CO(p,q) \ltimes \R^{p,q}$, and $i^+$ (respectively $i^-$) intertwines the adjoint action of $P$ on $\lieu^+$ (respectively $\lieu^-$) with the conformal action of $CO(p,q)$ on $\R^{p,q}$.


\subsubsection{Translations in $\mbox{PO}(p+1,q+1)$}
\label{section.translations}
Let $U^+$ be the closed subgroup of $\mbox{PO}(p+1,q+1)$ with Lie algebra $\lieu^+$.

\begin{defn}
A translation of $\mbox{PO}(p+1,q+1)$ is an element which is conjugate in $\mbox{PO}(p+1,q+1)$ to an element of $U^+$.
A translation of $\oo(p+1,q+1)$ is an element generating a $1$-parameter group of translations of $\mbox{PO}(p+1,q+1)$.
\end{defn}

This terminology is justified because a translation is a conformal transformation of $\Ein^{p,q}$ fixing a point, say $x$, and acting as a translation 
on ${\bf M}(x)$.  Notice that there are three conjugacy classes of translations in $\mbox{O}(p+1,q+1)$: lightlike (we will also say  null), spacelike, and timelike. 
An example of a null translation is the element $T=(i^+)^{-1}(1,0,\ldots,0)$ of $\lieu^+$.




Since any null translation of $\p$ is conjugate under $P$ to $T$, the reader will easily check the following fact, that  will be used several times below.
\begin{fact}
\label{centralizer-translation}
Let $T \in \p$ be a nontrivial null translation and $\mathfrak{c}(T)$ the centralizer  of $T$ in $\oo(p+1,q+1)$.  Then $\mathfrak{c}(T) \cap \p$ is of codimension  one in $\mathfrak{c}(T)$. 
\end{fact}

\subsection{Bounds in $\mbox{PO}(p+1,q+1)$}
\label{section.borne} The first step for proving theorem \ref{degree.bound} is to show that any nilpotent subalgebra of $\oo(p+1,q+1)$ has degree $\leq 2p+1$. We will actually prove more:

\begin{propn}
\label{opq.deg.bound}
For a nilpotent subalgebra $\h \subset \oo(p+1,q+1)$, the degree $d(\h) \leq 2p+1$. Assuming $p \geq
1$, if $d(\h) = d \geq 2p$, then $\h$ contains a translation in its center; in fact, $\h_{d-1}$ consists of null translations.
\end{propn}

\subsubsection{Preliminary results} The following definitions will be relevant below.  Let $\mathfrak{l} \subset \gl(n)$ be a subalgebra.  The set of all compositions $\Pi_1^k X_i$, where $X_1, \ldots, X_k \in \mathfrak{l}$, will be denoted $\mathfrak{l}^k$.  We say that $\mathfrak{l}$ is a {\it subalgebra of nilpotents} if there exists $k \geq 1$ such that $\mathfrak{l}^k = 0$. The minimal such $k$ will be called the \emph{order of nilpotence} of $\mathfrak{l}$, denoted $o(\mathfrak{l})$.  By Lie's theorem, subalgebras of nilpotents coincide with those subalgebras of  $\gl(n)$, the elements of which are nilpotent matrices.  If $\h$ is a nilpotent Lie algebra, then $\ad \h \subset \mathfrak{gl}(\h)$ is a subalgebra of nilpotents and $d(\h) = o(\ad \h)$.

For $V$ a vector space with form $B$, a Lie subalgebra $\h \subset \co(V)$ is \emph{infinitesimally conformal} if for all $u,v \in V$ and $X \in \h$,
$$ B( Xu,v) + B(u,Xv) = \lambda(X) B(u,v)$$ 
for some infinitesimal character $\lambda : \h \rightarrow \R$.  Of course, the Lie algebra of a subgroup of $\CO(V) \subset \GL(V)$ acts by infinitesimally conformal endomorphisms of $V$.

\begin{lmm}
\label{ad.of.brackets}
Let  $\mathfrak{l}$ be a Lie algebra and $V$ a finite dimensional  $\mathfrak{l}$-module. Let $Y \in \h_{k-1}$ and  $v \in V$.  Then $Y(v) \in {\mathfrak l}^k (V)$.
\end{lmm}

This lemma is easily proved by induction, using the Jacobi identity.



\begin{lmm}
\label{rep.order.bounds}
Let $V$ be a vector space with a symmetric bilinear form $B$ of type $(p,q)$. If $\overline{\lieu} \subset \co(V)$ is a subalgebra of nilpotents, then $o(\overline{\lieu}) \leq 2p+1$.
\end{lmm}

\begin{Pf}
Note that $\overline{\lieu}$ is infinitesimally isometric because there are no nontrivial infinitesimal characters $\overline{\lieu} \rightarrow \R$.
If $p=0$, then $\overline{\lieu} \subset \oo(n)$, in which case $\overline{\lieu}$ must be trivial.  

Now assume $p \geq 1$.  Let $U$ be the connected group of unipotent matrices in $CO(p,q)$ with Lie algebra $\overline{\lieu}$. Because $U$ consists of unipotent matrices, it lies
in a minimal parabolic subgroup of $\CO(p,q)$,  hence leaves invariant some isotropic $p$-plane $N \subset V$.
%
The order  of $\overline{\lieu}$ on both $N$ and $V/N^\perp$ is at most $p$, because each is dimension $p$.  Because $N^\perp/N$ inherits a positive-definite inner product that is infinitesimally conformally invariant by $\overline{\lieu}$, the order of $\overline{\lieu}$ on it is $1$.  Then $o(\overline{\lieu}) \leq 2p+1$, as desired.
\end{Pf}

\subsubsection{Proof of  proposition \ref{opq.deg.bound}}
\label{section.proof.deg.bound}
Let $H$ be the connected subgroup of $G = \mbox{PO}(p+1,q+1)$ with Lie algebra $\h$.  We recall  some facts from the theory of algebraic groups.  First, the nilpotence degrees of a connected group and its Zariski closure are the same, so that  there is no loss of generality assuming $H$ Zariski closed.  Then there is an algebraic Levi decomposition of the Lie algebra $\h \cong \mathfrak{r} \times \mathfrak{u}$, where $\mathfrak{r}$ is abelian and comprises the semisimple elements of $\h$, and  $\mathfrak{u}$ consists of nilpotents (see for example \cite[thm 10.6]{borel}).  If $\mathfrak{u}$ is trivial, $\h$ is abelian and proposition \ref{opq.deg.bound} is proved. If not, $d(\h) = d(\mathfrak{u})$, so we will assume that $\h=\mathfrak{u}$.  Moreover, $\mathfrak{u}$ is contained in a minimal
parabolic subalgebra of $\oo(p+1,q+1)$, and so conjugating $\mathfrak{u}$ if necessary, we have  $\lieu \subset \p$. Thanks to $i^+$ (see \ref{brief.description}), we identify $\mathfrak{u}$ with a subalgebra of $\co(p,q) \ltimes \R^{p,q}$.

Denote by $\overline{\mathfrak{u}}$ the projection to $\co(p,q)$, which is actually in $\oo(p,q)$ since $\overline{\mathfrak{u}}$ is a subalgebra of nilpotents.  For any natural number $k$,
\begin{equation}
\label{eq.nilpotence}
 \mathfrak{u}_k \subseteq \overline{\mathfrak{u}}_k + \overline{\mathfrak{u}}^k (\R^{p,q})
 \end{equation}
 
The proof by induction of this relation is straightforward using lemma \ref{ad.of.brackets} for $\mathfrak{l}=\mathfrak{u}$ and $V=\R^{p,q}$, and is left to the reader.  

When $p=0$, then any nilpotent subalgebra $\overline{\mathfrak{u}} \subset \oo(1,q+1)$ is abelian, by the remarks above combined with lemma \ref{rep.order.bounds}. 
We thus have $d(\lieu) \leq 2p+1$ when $p=0$.  Now proceed inductively on $p$, using lemma \ref{rep.order.bounds} and relation (\ref{eq.nilpotence}) to obtain $d(\mathfrak{u}) \leq o(\overline{u}) \leq 2p+1$ whenever $\lieu \subset \oo(p+1,q+1)$ is nilpotent, for all $p \in \N$.

Next suppose that $ d(\mathfrak{u}) = d \geq 2p \geq 2$.  Since $\overline{\mathfrak{u}}$ is a nilpotent subalgebra of $\oo(p,q)$, its nilpotence degree is at most $2p-1$ by the first part of the proof.  Since $d \geq 2p$, ${\overline{\mathfrak{u}}}_{d-1} = 0$ and 
$$ 0 \neq \mathfrak{u}_{d-1} \subseteq \overline{\mathfrak{u}}^{d-1} (\R^{p,q})$$
so an element of $\mathfrak{u}_{d-1}$ can be written 
$$w = Y_1 \cdots Y_{d-1} (v) \qquad \mbox{for} \ Y_1, \ldots, Y_{d-1} \in \overline{\mathfrak{u}}, \ v \in \R^{p,q}$$

Further, any $Y \in \overline{\mathfrak{u}}$ annihilates $w$.  Because $Y_1 \in \overline{\mathfrak{u}}$ is infinitesimally conformal and nilpotent, it is infinitesimally isometric.  Then
$$ Q^{p,q}(w,w) = \langle Y_1 \cdots Y_{d-1} (v), Y_1 \cdots Y_{d-1} (v) \rangle = - \langle Y_2 \cdots Y_{d-1} (v), Y_1 (w) \rangle = 0$$
and so $w$ is a null translation.
$\diamondsuit$

\subsection{Conformal structures as Cartan geometries}
\label{section.cartan}
In the sequel, it will be fruitful to study pseudo-Riemannian structures in the setting of Cartan geometries.  A Cartan geometry modeled on some homogeneous space $\X=G/P$ is a curved analogue of $\X$. 
\begin{defn}
\label{definition.cartan} Let $G$ be a Lie group with Lie algebra $\g$ and $P$ a closed subgroup of $G$ such that $\Ad P$ is faithful on $\g$.  A \emph{Cartan geometry} $(M,B,\omega)$ modeled on $(\g,P)$ is
\begin{enumerate}
\item{a principal $P$-bundle $\pi: B \rightarrow M$}
\item{a $\g$-valued $1$-form $\omega$ on $B$ satisfying}
\begin{itemize}
\item{for all $b \in B$, the restriction $\omega_b : T_bB \rightarrow \g$ is an isomorphism} 
\item{for all $b \in B$ and $Y \in \h$, the evaluation $\omega_b(\left. \frac{\mbox{d}}{\mbox{d} t} \right|_0 be^{tY} ) = Y$ }
\item{for all $b \in B$ and $h \in P$, the pullback $R_h^* \omega = \Ad h^{-1} \circ \omega$}
\end{itemize}
\end{enumerate}
\end{defn}

For the model $\X$, the canonical Cartan geometry is the triple $(\X,G,\omega_G)$, where $\omega_G$ denotes the left-invariant $\g$-valued $1$-form on $G$, called the \emph{Maurer-Cartan form}.

A conformal structure $(M,[\sigma])$ of type $(p,q)$ with  $p+q \geq 3$,  defines, up to isomorphism, a canonical Cartan geometry  $(M,B,\omega)$ modeled on $\Ein^{p,q}$---that is, on $(\oo(p+1,q+1),P)$.  The interested reader will find the details of this solution, originally due to E. Cartan, of the so-called {\it equivalence problem} in \cite[ch 7]{sharpe}.  

The group $\Aut M$ comprises the bundle automorphisms of $B$ preserving $\omega$.  Any conformal diffeomorphism lifts to an element of $\Aut M$, maybe not unique, but the fibers of the projection from  $\Aut M$ to the conformal group of $M$ are discrete.  The Lie algebras of $\Aut M$ and $\Conf M$ are thus isomorphic.  We will not distinguish in notation between an element $f \in \Conf M$ and the corresponding lift to $\Aut M$.    

\section{General degree bound: proof of theorem \ref{degree.bound}}
\label{section.general.degree.bound}
In this section, we  use  the interpretation of conformal structures as Cartan
geometries to prove theorem \ref{degree.bound}.  Let $(M,B,\omega)$ be a Cartan geometry modeled on $(\g, P)$, and $H < \Aut M$ a
connected Lie group.  Since $H$ acts on $B$, each vector $X \in \h$ defines a Killing field on $B$, and for every $b \in B$, we will call $X(b)$ the value of this Killing field at $b$.  Thus,  each point $b \in B$ determines a linear embedding
\begin{eqnarray*}
 s_b & : & \h \to \g \\
& &  X \mapsto \omega_b(X(b))
\end{eqnarray*}

The injectivity of $s_b$ comes from the fact that $H$ preserves a framing on
$B$, hence acts freely (see \cite[I.3.2]{kob}).  The image $s_b(\h)$ will be denoted $\h^b$, and, for
$X \in \h$, the image $s_b(X)$ will be denoted $X^b$.  In
general, $s_b$ is not a Lie algebra homomorphism, except with respect to stabilizers (see \cite[5.3.10]{sharpe}): for any $X,Y \in \h$ and $b \in B$ such that $Y^b \in \p$, 
$$ [X,Y]^b = [X^b,Y^b]$$
  Observing that $Y$ belongs to the stabilizer $\h(\pi(b))$ if
and only if $Y^b \in \p$, we deduce the following fact.
\begin{fact}
\label{codim-one}
If $\h^b \cap \p$ is codimension at most $1$ in $\h^b$, then $\h^b$ is a Lie subalgebra of $\g$, isomorphic to $\h$.
\end{fact}

The following result implies theorem \ref{degree.bound}. It is more precise and will be useful for the proof of theorem \ref{flatness.thm}:
\begin{thm}
\label{degree.bound.precise}
Let $(M,[\sigma])$ be a compact manifold with a type-$(p,q)$ conformal structure, and let
$(M,B,\omega)$ be the associated Cartan geometry.  Let $H <
\Aut M$ be a connected nilpotent Lie group.  Then $d(H) \leq 2p+1$.  If $d(H)
= 2p+1$, then every $H$-invariant closed subset  $F \subset M$ contains a
point $x$ such that
\begin{enumerate}
\item  The dimension of the orbit $H.x$ is at most $1$.

\item For every $b \in \pi^{-1}(x)$, $\h^b$ is a subalgebra of $\oo(p+1,q+1)$.

\item There exists $X \in \h$ such that $X^b$ is a lightlike translation in
  $\p$ for every $b \in \pi^{-1}(x)$, and $X^b$ is in the center of $\h^b$.
\end{enumerate}
\end{thm}

A consequence of this theorem is that when $d(H) = 2p+1$, there are points
with nontrivial stabilizers, because $X$ as in (3) generates a $1$-parameter
subgroup of the stabilizer $H(x)$.  We will study the dynamics near $x$ of this flow in the proof of theorem \ref{flatness.thm}.




\begin{Pf} (of theorem \ref{degree.bound.precise})

Let $F \subset M$ be closed and $H$-invariant.  The group $H$ is amenable, so it preserves a finite Borel measure on $F$.  Then the embedding theorem of \cite[thm 4.1]{bfm} with $S = H$ gives $x \in F$ and an algebraic subgroup $\check{S} < \mbox{Ad}_\g P$ such that, for all $b \in \pi^{-1}(x)$,                              
\begin{enumerate}
\item $\h^b$ is $\check{S}$-invariant
\item $s_b$ intertwines the Zariski closure of $\Ad H$ in $\Aut \h$ with $\left. \check{S} \right|_{\h^b}$ 
\end{enumerate}

Since the adjoint representation of $\mbox{PO}(p+1,q+1)$ is
algebraic and faithful, $\check{S}$ is the image of an algebraic subgroup of
$P$, which we will also denote $\check{S}$.  We denote the corresponding Lie algebra by $\check{\mathfrak{s}} \subset \mathfrak{p}$.  The embedding theorem says that for any $X \in \h$, there exists $\check{X} \in
\check{\mathfrak{s}}$ such that for all $Y \in \h$,
$$ [X,Y]^b = [\check{X},Y^b]$$

Suppose that $d = d(\h) \geq 2p+1$.
Because $\check{\s}$ is algebraic, there is a decomposition $\check{\s} \cong \mathfrak{r} \ltimes \mathfrak{u}$ 
with $\mathfrak{r}$ reductive and $\mathfrak{u}$ consisting of nilpotent
elements (see \cite[4.4.7]{wm}).  Because $\ad \h$ consists of nilpotents, the subalgebra
$\mathfrak{r}$ is in the kernel of restriction to $\h^b$, and $\lieu$ maps
onto $\ad \h$.  Therefore, for $l = d(\lieu)$, 
$$ 2p \leq d-1 = d(\ad \h) \leq l \leq 2p+1$$ 
where the upper bound comes from proposition \ref{opq.deg.bound}. Also by this proposition, $\lieu_{l-1}$ consists of null translations.  Whether $l=d-1$ or $d$, we will show that $\h^b$ centralizes a null translation in $\p$, from which fact we will obtain the bound and points (1) and (3).

First suppose $l=d-1$.  Then $\lieu_{d-2}$ consists of null translations and acts on $\h^b$ as $\ad \h_{d-2}$, which means it centralizes $(\h_1)^b$.  Then by facts \ref{centralizer-translation} and \ref{codim-one}, $(\h_1)^b$ embeds homomorphically in $\mathfrak{o}(p+1,q+1)$.  The order of $\mathfrak{u}$ on $(\h_1)^b$ is $d-1$; further, $\mathfrak{u}$
and $(\h_1)^b$ generate a nilpotent subalgebra $\mathfrak{n}$ of order $d-1$, in which $(\h_1)^b$ is an ideal.  Since $d-1 \geq 2p$,
proposition \ref{opq.deg.bound} implies that the commutators
$\mathfrak{n}_{d-2}$ are all null translations.  But $\lien_{d-2}$ contains
$$ \mathfrak{u}^{d-2} (\h_1)^b = (\h_{d-1})^b $$
Because $\mathfrak{u}$ preserves $(\h_1)^b \cap \p$ and acts by nilpotent transformations on $(\h_1)^b / ((\h_1)^b \cap \p)$, which is $1$-dimensional,  
$$\mathfrak{u}^1 (\h_1)^b = (\h_2)^b \subset \p$$
Thus $(\h_k)^b \subset \p$ as soon as $k \geq 2$, so for any
$X \in \h$, $Y \in \h_k$, we have $[X,Y]^b=[X^b,Y^b]$.  In particular,
$(\h_{d-1})^b$, the image under $s_b$ of the center of $\h$, commutes with $\h^b$, so that $\h^b$ is in the centralizer of a
nonzero null translation.

Next suppose $l=d$.  Then $\lieu_{d-1}$ centralizes $\h^b$ because it acts as $\ad \h_{d-1}$.  By proposition \ref{opq.deg.bound}, $\lieu_{d-1}$ consists of null translations, so $\h^b$ commutes with a nonzero null translation in $\p$.

Given that $\h^b$ centralizes a nonzero null translation in $\p$, \ref{centralizer-translation} implies point (1); moreover, $(\h^b)_1 \subset \p$.  By fact \ref{codim-one}, $s_b : \h \rightarrow \oo(p+1,q+1)$ is a homomorphic embedding.  The assumption $d \geq 2p+1$ and proposition \ref{opq.deg.bound} forces $d = 2p+1$, proving the bound.  Also by proposition \ref{opq.deg.bound}, $(\h^b)_{2p}$, which is central in $\h^b$, consists of null translations.  Then $(\h^b)_{2p} \subset (\h^b)_1 \subset \p$, and (3) is proved.

Finally, item $(2)$ of the theorem is a consequence of item (1) and fact \ref{codim-one}. \end{Pf}

\subsection{Outline of the proof of theorem \ref{flatness.thm}}

For $x \in M$, denote by $H(x)$ the stabilizer of $x$ in $H$.  
For each $b \in \pi^{-1}(x)$, the action of $H$ by automorphisms of the
principal bundle $B$ gives rise to an injective homomorphism $\rho_b : H(x)
\rightarrow P$.  Theorem \ref{degree.bound.precise} says that if $d(H) = 2p+1$, then each
$H$-invariant closed set $F$ contains a point $x_0$, such that for some
$1$-parameter group $h^s$ in $H(x_0)$ and $b_0 \in \pi^{-1}(x_0)$, the image $\rho_{b_0}(h^s)$ is a
$1$-parameter group $\tau^s$ of null translations in $P$.  

The dynamics of null translations are studied in section \ref{subsection.dynamics.null} and summarized in fact \ref{fact.taus.dynamics}.  In section \ref{geodesics.holonomy}, we make
the crucial link between the dynamics of $\tau^s$ on $\Ein^{p,q}$ and those of
$h^s$ on $M$, via the respective actions on special curves, called \emph{geodesics}, in the two Cartan
bundles.  The actions on these curves are conjugate
locally by the \emph{exponential maps} of the two Cartan geometries.  In section \ref{subsection.dynamics.on.M}, we deduce from this relationship the dynamics of the $h^s$-action and develop a method to precisely compute the differential of $h^s$ near $x_0$; see proposition \ref{prop.combined.reparam.framing}. 

In section \ref{vanishing.curvature} we use our description of the $h^s$-action to show that $M$ is {\it conformally flat}, namely locally modeled on $\Ein^{p,q}$ (see proposition \ref{prop.step.flatness}).  First we reduce the claim to showing flatness on a neighborhood of $x_0$.  Proposition \ref{prop.cartan.vanishing} of section \ref{sec.vanishing.lightcone} establishes flatness on a nonempty open subset.  In section \ref{sec.vanishing.neighborhood}, we make a technical modification of proposition \ref{prop.cartan.vanishing} in order to show that this flat set includes a neighborhood of $x_0$ (see proposition \ref{prop.vanishing.on.nbhd}).

The purpose of section \ref{global} is to understand the global structure of $M$, using classical techniques for $(G,X)$-structures.  We again use the dynamics of $h^s$ to show that the subset of $M$ tending under the forward or reverse flow to a fixed curve develops to the complement of the fixed set of $\tau^s$ in $\Ein^{p,q}$; the Lorentzian case is treated first in section \ref{sec.proof.lorentzian}, while the more complicated case $p \geq 2$ is in section \ref{sec.proof.p}.  In both cases, we apply the the theorem \cite[thm 1.8]{charles.bords} on boundaries of embeddings of flat Cartan geometries to conclude.  This final section completes the proof of theorem \ref{flatness.thm}.

\section{Conformal dynamics}
\label{conformal.dynamics}

This section establishes properties of the flow $h^s$ that will later be used to show $M$ is conformally flat.  Along the way, we develop some general tools to relate the behavior of an automorphism of a Cartan geometry with the behavior of the corresponding holonomy on the model space.


\subsection{Dynamics of null translations on $\Ein^{p,q}$}
\label{subsection.dynamics.null}
The first task is to describe the dynamics of $1$-parameter groups of null translations in the model space $\Ein^{p,q}$.  Let $\tau^s$ be the flow generated by the null translation $T = (i^+)^{-1}(1, 0 , \ldots, 0)$ of section \ref{section.translations}.


The action of $\tau^s$ on $\Ein^{p,q}$ is given in projective coordinates by
$$ \tau^s  :  \ [ y_0 : \  \cdots \ : y_{n+1}] \mapsto [y_0 + sy_n : y_1 - sy_{n+1}: y_2 : \ \cdots \ : y_{n+1}] $$

The fixed set is 
$$F= \BP(e_0^\perp \cap e_1^\perp \cap \mathcal{N}^{p+1,q+1})$$

When $p=2$, it has codimension $2$, and contains a singular circle 
$$\Lambda= \BP(\mbox{span} \{ e_0,e_1 \}) \subset F$$
When $p=1$, then $F = \Lambda$; since $p+q \geq 3$, the codimension of $F$ is also at least $2$ in this case.

If $y \notin F$, then 
$$\tau^s. y \rightarrow [y_n: -y_{n+1} : 0 : \cdots : 0] \in \Lambda \qquad \mbox{as}\ s \rightarrow \infty$$  

\begin{figure}[ht]
\includegraphics[width=3in,height=2in]{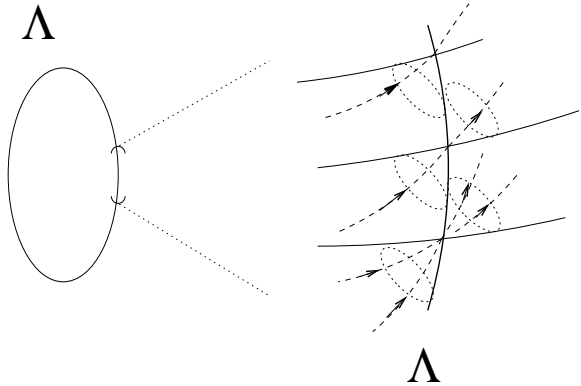}
\caption{local picture of flow by null translation $\tau^s$}
\protect{\label{figure.nullflow}}
\end{figure}

Every point $x \in \Ein^{p,q}$ lies in some $C(y)$ for $y \in \Lambda$, and
$y$ is unique when $x \notin F$.  We summarize the dynamics of $\tau^s$ near
$\Lambda$; see also figure ~\ref{figure.nullflow}:

\begin{fact} 
\label{fact.taus.dynamics}
The complement of the closed, codimension-$2$ fixed set $F$ of
  $\tau^s$ in $\Ein^{p,q}$ is foliated by subsets of lightcones $\check{C}(y) = C(y) \setminus (C(y) \cap
  F)$, for $y \in \Lambda$. Points $x \in \check{C}(y)$ tend under $\tau^s$ to $y$
  along the lightlike geodesic containing $x$ and $y$; in particular, $\tau^s$
preserves setwise all null geodesics emanating from points of $\Lambda$.
\end{fact}

\subsection{Geodesics and holonomy}
\label{geodesics.holonomy}



In this section $(M,B,\omega)$ will be a Cartan geometry modeled on $G/P$.  The form $\omega$ on $B$ determines special curves, the \emph{geodesics}.  Here they will be defined as projections of curves with constant velocity according to $\omega$---that is, $\gamma : (- \epsilon, \epsilon) \rightarrow M$ is a geodesic if $\gamma(t) = \pi(\hat{\gamma}(t))$ where
$$ \omega(\hat{\gamma}'(t)) = \omega(\hat{\gamma}'(0)) \qquad \mbox{for all}\ t \in ( -\epsilon, \epsilon)$$

Geodesics on the flat model space $(G/P,G,\omega_G)$ are orbits of
$1$-parameter subgroups.  Note that this class of curves is larger than the
usual set of geodesics in case the Cartan geometry corresponds to a pseudo-Riemannian
metric or a conformal pseudo-Riemannian structure (see \cite{fialkow}, \cite{friedrich}, \cite{friedrich.schmidt} for a definition of
conformal geodesics).

The \emph{exponential map} is defined on $B \times \g$ in a neighborhood of $B \times \{ 0 \}$ by
$$ \mbox{exp}(b,X) = \mbox{exp}_b(X) = \hat{\gamma}_{X,b}(1)$$

where $\hat{\gamma}_{X,b}(0) = b$ and $\omega(\hat{\gamma}_{X,b}'(t)) = X$ for
all $t$.

Let $h \in \aut M$.  Then $h$ carries geodesics in $M$ to geodesics:  

$$\hat{h} \circ \hat{\gamma}_{X,b} = \hat{\gamma}_{X,\hat{h}(b)}$$

Suppose that $h^s$ is a $1$-parameter group of
automorphisms 
Then for any $b_0 \in B$, the curve parametrized by the flow
$\hat{\gamma}(s) = h^s.b_0$ projects to a
geodesic $\gamma(s)$ in $M$: it is easy to see that $\hat{\gamma}(s)$ has $\omega$-constant velocity.


\begin{defn}
If $h \in \aut M$ fixes $x$, and $b \in \pi^{-1}(x)$, then the element $g \in P$ such that $h.b=bg$ is the \emph{holonomy of $h$ with respect to $b$}.   More generally, given a local section $\sigma : U \rightarrow B$, the \emph{holonomy of $h \in \aut M$ with respect to $\sigma$} at some point $x \in U \cap h^{-1}.U$ is $g$ such that $h. \sigma(x) = \sigma(h.x) g$ 
\end{defn}


If $H(x)$ is the stabilizer of $x$ in $H < \aut M$ and $b \in \pi^{-1}(x)$, then the holonomy with respect to $b$ gives a monomorphism
$\rho_b : H(x) \rightarrow P$.  Replacing $b$ with $bp$ has the effect of
post-composing with conjugation by $p^{-1}$. 

For automorphisms fixing a point $x_0$, the holonomy with respect to some $b_0 \in
\pi^{-1}(x_0)$ tells a lot about the action in a neighborhood of $x_0$ via the
exponential map.  If, moreover, an automorphism $h$ fixes $x_0$ and preserves
the image of a geodesic $\gamma$ emanating from $x_0$, then the holonomy at $x_0$ determines the holonomy along
$\gamma$, as follows.  For $X \in \g$, denote by $e^X$ the exponential of $X$ in $G$.

\begin{propn}
\label{prop.holonomy}
Suppose that $h \in \aut M$ fixes a point $x_0$ and has holonomy $g_0$ with respect
to $b_0 \in \pi^{-1}(x_0)$.  Let $\gamma(t) = \pi(\exp(b_0,tX))$ for $X \in \g$,
defined on an interval $(\alpha, \beta)$ containing $0$.  Suppose there exist 
\begin{itemize}
\item a path $g : (\alpha,
\beta) \rightarrow P$ with $g(0) = g_0$
\item a diffeomorphism $c : (\alpha,
\beta) \rightarrow (\alpha', \beta')$
\end{itemize}
 such that, for all $t \in (\alpha, \beta)$,
$$ g_0 e^{tX} = e^{c(t)X} g(t) $$
Then
\begin{enumerate}
\item  The curves $\exp(b_0,c(t)X)$ and $\gamma(c(t))$ are defined for all $t \in (\alpha,
  \beta)$, and $h.\gamma(t) = \gamma(c(t))$.



\item Viewing $\exp(b_0,tX)$ as a
  section of $B$ over $\gamma(t)$, the holonomy of $h$ at $\gamma(t)$ with
  respect to this section is $g(t)$.



\end{enumerate}
\end{propn}

\begin{Pf}
In $G$, reading the derivative of $g_0 e^{tX}$ with $\omega_G$ gives (see {\cite[3.4.12]{sharpe}})
\begin{eqnarray*}
X = (\Ad g(t)^{-1})(c'(t)X) + \omega_G(g'(t))
\end{eqnarray*}

Because $h$ preserves $\omega$, the derivative of 
$$h.\exp(b_0,tX) = \exp(b_0 g_0, tX)$$ 

according to $\omega$ is $X$ for all $t \in (\alpha, \beta)$.  On the other hand, it is also true in $B$ that whenever $t \in (\alpha, \beta)$ and $\exp(b_0,c(t)X)$ is defined,
$$ \omega((\exp(b_0,c(t)X) g(t))') =  (\Ad g(t)^{-1})(c'(t)X) + \omega_G(g'(t))$$

This formula follows from the properties of $\omega$ in the definition \ref{definition.cartan} of a Cartan geometry; see {\cite[5.4.12]{sharpe}}.  Therefore, because the two curves have the same initial value and the same derivatives, both are defined for $t \in (\alpha,\beta)$, and
$$ h.\exp(b_0,tX) = \exp(b_0,c(t)X) g(t)$$
This proves (2); item (1) follows by projecting both curves to $M$. \end{Pf}



We record one more completeness result that will be useful below, for flows
that preserve a geodesic, but do not necessarily fix a point on it.  


\begin{propn}
\label{prop.flow.completeness}
Let $X \in \g$ and $b_0 \in B$ be such that $\exp(b_0,tX)$ is defined for all $t \in (\alpha,\beta)$, where $\alpha < 0 < \beta$.
Suppose that for some $Y \in \g$, there exists $g: (\alpha,\beta) \rightarrow P$ such that for all $t \in (\alpha, \beta)$,
$$ e^{tX} = e^{c(t)Y}  g(t)$$
in $G$, where $c$ is a diffeomorphism $(\alpha,\beta) \rightarrow \R$ fixing $0$.  Then $\exp(b_0,tY)$ is defined for all $t \in \R$.
\end{propn}

\begin{Pf}
In $G$, we have for all $t \in (\alpha,\beta)$
\begin{eqnarray}
\label{eqn.y.equals.derivative}
Y = \frac{1}{c'(t)} [\Ad g(t)][X - \omega_G(g'(t))]
\end{eqnarray}
 Let $c^{-1}(s)$ be the inverse diffeomorphism $(- \infty,\infty) \rightarrow
 (\alpha, \beta)$.
Define, for $s \in \R$,
$$ \hat{\gamma}(s) = \exp(b_0,c^{-1}(s)X)  g(c^{-1}(s))^{-1}$$
The derivative of the right-hand side is
$$\frac{1}{c'(c^{-1}(s))}  [\Ad g(c^{-1}(s))] [X - \omega_G(g'(c^{-1}(s))]$$
For $t = c^{-1}(s)$, this reduces to the right-hand side of (\ref{eqn.y.equals.derivative}).
Then $\exp(b_0,sY)$ equals $\hat{\gamma}(s)$ and is defined for all $s \in \R$. \end{Pf}






\subsection{Dynamics of  $h^s$ on $M$}
\label{subsection.dynamics.on.M}

We now return to the pseudo-Riemannian manifold $(M,\sigma)$ with associated Cartan geometry $(M,B, \omega)$ modeled on the pair $(\oo(p+1,q+1),P)$ corresponding to the homogeneous space $\Ein^{p,q}$.  The subalgebra $\lieu^-$ complementary to $\p$ and the basis $U_1, \ldots, U_n$ are as in 
section \ref{brief.description}. Let $\mathcal{N}(\mathfrak{u}^-)$ be the null cone with respect to $Q^-:=(i^-)^*(Q^{p,q})$ in $\mathfrak{u}^-$ (see \ref{brief.description} for the definition of $i^-$).  The following proposition captures properties of the flow $h^s$ given by theorem \ref{degree.bound.precise} that reflect properties of $\tau^s$ established above, and that will be used to prove flatness in the next section.

Now we suppose the group $H < \Conf M$ is nilpotent of maximal degree $2p+1$.  The point $x_0$ is given by theorem \ref{degree.bound.precise}, and $h^s.x_0=x_0$.
We have $b_0 \in \pi^{-1}(x_0)$ for which the holonomy of $h^s$ is $\tau^s$ as in section \ref{subsection.dynamics.null}.


\begin{propn}
\label{prop.combined.reparam.framing}
Suppose $H < \Conf M$ is nilpotent of maximal degree $2p+1$.  Let $F \subseteq M$ be closed and $H$-invariant subset, and let $x_0 \in F$ and $X \in \h$ be given by theorem \ref{degree.bound.precise}.  Let $h^s$ be the flow generated by $X$. Then there is $b_0 \in \pi^{-1}(x_0)$ for which the following holds: 
\begin{enumerate}
\item  Let $\hat{\Delta}(v)
= \mbox{exp}(b_0,v U_1)$ with domain $I_\Delta \subseteq \R$; let $\Delta = \pi \circ \hat{\Delta}.$  Then the flow $h^s$ fixes pointwise the geodesic $\Delta$, and, for $v \in I_{\Delta}$, has holonomy at $\Delta(v)$ with respect to $\hat{\Delta}$ equal ${\tau}^s$.

\item  There is an open subset $\mathcal{S} \subset \mathcal{N}(\lieu^-)$ such
that $\mathcal{S} \cup - \mathcal{S}$ is dense in $\mathcal{N}(\lieu^-)$, and for all $v \in I_\Delta$ and $U \in \mathcal{S}$,
if the geodesic $\beta(t) = \pi \circ
    \exp(\hat{\Delta}(v), t U)$ is defined on $(-\epsilon,\epsilon)$, then the flow $h^s$ preserves $\beta$ and reparametrizes by
$$ c(t) = \frac{t}{1+st}$$
for $t \in (-\epsilon,\epsilon)$.  In particular, for $t > 0$ \ ($t <0$), 
$$h^s(\beta(t)) \rightarrow \Delta(v) \ \mbox{as} \ s \rightarrow \infty \ (s \rightarrow - \infty) $$
Moreover, $\beta(t)$ is complete.

\item There is a framing $f_1(t), \ldots, f_n(t)$ of $M$ along
$\beta(t)$ for which the
derivative
$$ h^s_*(f_i(t)) = \left( \frac{1}{1+st}\right) ^{\sigma(i)} f_i(c(t))$$
where 
$$ \sigma(i) = \left\{
\begin{array}{ccl}
0 & \qquad & i = 1 \\
1 & \qquad & i \in \{ 2, \ldots, n-1 \} \\
2 & \qquad & i = n
\end{array} \right.
$$

\end{enumerate}
\end{propn}

To prove this proposition, we will compute the relevant holonomies of $\tau^s$ in the model space and use proposition \ref{prop.holonomy}. 
We start with some algebraic facts pertaining to $\Ein^{p,q}$.  Let $\lier$ be as in section \ref{brief.description}, a
maximal reductive subalgebra of $\p$.
\begin{lmm}
\label{lmm.alg.facts}
Let $R \cong CO(p,q)$ be the connected subgroup of $P$ with Lie algebra
$\mathfrak{r}$, and let $S$ be the unipotent radical of the stabilizer in $R$ of
$U_1$.  
\begin{enumerate}
\item{$\mbox{Fix}(\Ad \tau^s) \cap \mathfrak{u}^- = \R U_1$}
\item{Let $\mathcal{S} = \R_{>0}^* \cdot S.U_n$.  Then $\mathcal{S} \cup - \mathcal{S}$ is open and dense in $\mathcal{N}(\lieu^-)$.}
\item{The subgroups $S$ and $\tau^s$ commute.}
\end{enumerate}
\end{lmm}

\begin{Pf}
\begin{enumerate}
\item Recall that $T$ is the infinitesimal generator for $\tau^s$ defined in section \ref{section.translations}.  It
  suffices to show 
$$ \R U_1 = \ker(\ad T) \cap \lieu^-$$
We leave this basic linear algebra calculation to the reader.

\item We will show that $\mathcal{S}$ consists of all $U \in
  \mathcal{N}(\lieu^-)$ with $\langle U, U_1 \rangle > 0$; these elements and their negatives form an
  open dense subset of $\mathcal{N}(\lieu^-)$. 

First, if $g \in S$, then $g. U_n \in \mathcal{N}(\lieu^-)$ and
\begin{eqnarray*}
\langle g. U_n ,U_1 \rangle & = & \langle g. U_n, gU_1 \rangle =  \langle U_n,U_1 \rangle = 1
\end{eqnarray*}
Both $\mathcal{N}(\lieu^-)$ and the property $\langle U,U_1 \rangle > 0$ are invariant by the action of $\R^*_{>0}$, so $\mathcal{S}$ is contained in the claimed subset.

Next let $U \in \mathcal{N}(\lieu^-)$ be such that $\langle U,U_1 \rangle > 0$.  Replace $U$ with a positive scalar multiple so that $\left< U, U_1 \right> =1$.  Define $g \in S$ by
\begin{eqnarray*}
g & : & U_1 \mapsto U_1 \\
 &  &  U_n \mapsto U \\
 &  & V \mapsto V - \left<V,U \right>\cdot U_1 \qquad \mbox{for} \ V \in \{ U_1, U_n
 \}^\perp
\end{eqnarray*}
It is easy to see that $g$ is unipotent and belongs to $O(Q^-)$, and thus defines an element of $R$.  Therefore $U \in \mathcal{S}$.

\item Both $S$ and $\tau^s$ lie in the unipotent radical of $P$, which, in the
  chosen basis, is contained in the group of upper-triangular matrices.  The
  commutator of any unipotent element with $\tau^s$ is $I_{n+2} + c E_0^{n+1}$
  for some $c \in \R$.  There is no such element of $\mbox{O}(p+1,q+1)$ for any
  nonzero $c$, so the commutator is the identity.  
\end{enumerate}
\end{Pf}

\begin{Pf} (of proposition \ref{prop.combined.reparam.framing})

Because $\Ad \tau^s$ fixes $U_1$, the corresponding $1$-parameter subgroups
commute in $G$:
$$ \tau^s e^{vU_1} = e^{vU_1} \tau^s$$
Then by proposition \ref{prop.holonomy}, the flow
$h^s$ fixes $\Delta$ pointwise, and the holonomy of $h^s$ with
respect to $\hat{\Delta}$ at any $\Delta(v), \ v \in I_\Delta$, equals $\tau^s$.  This proves (1).

To prove (2), first consider the null geodesic $\alpha(t) = \pi (e^{t U_n})$ in $G/P$, and let $\hat{\alpha}(t) = e^{tU_n}$.  Now it is possible to compute the holonomy of $\tau^s$ along $\alpha$ with respect to $\hat{\alpha}$:



\begin{eqnarray*}
\tau^s \cdot \hat{\alpha}(t) & = & \tau^s \cdot e^{tU_n} \\
 & = & \hat{\alpha}(c(t)) \cdot e^{-c(t) U_n} \cdot \tau^s \cdot e^{t U_n}
\end{eqnarray*}

Refer to the expression for $\tau^s$ in section \ref{subsection.dynamics.null}, and compute directly in $\mbox{O}(p+1,q+1)$ 
$$ e^{-c(t) U_n} \cdot \tau^s \cdot e^{t U_n} = \mbox{diag}(1+st,1+st, 1, \ldots, 1, \frac{1}{1+st}, \frac{1}{1+st}) + sT$$

Denote this holonomy matrix by $h(s,t)$.

Now let $S < G$ be as in lemma \ref{lmm.alg.facts}, and let $U = (\Ad g)(U_n)$ with $g \in S$.  Let $\hat{\alpha}(t) = e^{tU}$.  Because $\tau^s$ commutes
with $g$ by lemma \ref{lmm.alg.facts} (3), we can compute the holonomy of $\tau^s$ with respect to $\hat{\alpha}$ along $\alpha$:
\begin{eqnarray*}
\tau^s \cdot \hat{\alpha}(t) & = & \tau^s \cdot e^{t U}  =  \tau^s \cdot e^{(\Ad g)(t U_n)} \\
  & = & \tau^s \cdot g \cdot e^{t U_n} \cdot g^{-1}  =  g \cdot \tau^s \cdot e^{tU_n} \cdot g^{-1} \\
  & = & g \cdot e^{c(t) U_n} \cdot h(s,t) \cdot g^{-1} =  e^{(\Ad g)(c(t) U_n)} \cdot g \cdot h(s,t) \cdot g^{-1} \\
  & = & \hat{\alpha}(c(t)) \cdot g \cdot h(s,t) \cdot g^{-1}
\end{eqnarray*}

Let $\mathcal{S}$ be as in lemma \ref{lmm.alg.facts} (2).  Let $U \in \mathcal{S}$.  Let $\hat{\beta}(t) =
\exp(\hat{\Delta}(v),tU)$ and $\beta = \pi \circ \hat{\beta}$, and assume $\hat{\beta}$ is defined on $(-\epsilon,\epsilon)$.  From (1), the holonomy of $h^s$ at
$\Delta(v)$ with respect to $\hat{\Delta}$ is $\tau^s$.  The above calculation, together with
proposition \ref{prop.holonomy} (1), implies
$$h^s.\beta(t) = \beta(c(t))$$ 

for all $t \in (- \epsilon,\epsilon)$.  Taking $s = \pm 1/\epsilon$ and again applying proposition \ref{prop.holonomy} (1)
 proves completeness of $\beta(t)$.  Then point (2) is proved.

By proposition \ref{prop.holonomy} (2), the holonomy of $h^s$
at $\beta(t)$ with respect to $\hat{\beta}$ is $g \cdot h(s,t) \cdot g^{-1}$. 
The adjoint of $h(s,t)$ on $\g/\p$ in the basis comprising the images of $U_1, \ldots, U_n$  is 
$$ \mbox{diag}(1, \frac{1}{1+st}, \ldots, \frac{1}{1+st}, \frac{1}{(1+st)^2} ) $$

Since $S$ is contained in $P$, for $g \in S$, the span of $(\Ad g)(U_1) , \ldots, (\Ad g)(U_n)$ is transverse to $\p$.
The adjoint of $g \cdot h(s,t) \cdot g^{-1}$ in the corresponding basis of $\g/\p$ is of course the same diagonal matrix as for $g=1$.  For $\hat{\beta}$ and $\beta$ as above, define a framing $f_1, \ldots, f_n$ along $\beta$ by
$$ f_i(\beta(t)) = (\pi_* \circ \omega^{-1}_{\hat{\beta}(t)} \circ \Ad g)(U_i)$$

Now we can compute the derivative of $h^s$ along $\beta$ in the framing $(f_1,
\ldots, f_n)$.  Recall the identity for a Cartan connection
$$ \omega_p^{-1} \circ (\Ad g) = R_{g^{-1}*} \circ \omega_{pg}^{-1} $$

We will write $f_i(t)$ in place of $f_i(\beta(t))$ below.  
\begin{eqnarray*}
h^s_* (f_i(t)) & = & \left( \pi_* \circ h^s_* \circ
\omega_{\hat{\beta}(t)}^{-1} \circ \Ad g \right)(U_i) \\
& = & \left( \pi_* \circ \omega_{h^s \cdot \hat{\beta}(t)}^{-1} \circ \Ad g \right)(U_i) \\
& = & \left( \pi_* \circ R_{g^{-1}*} \circ \omega_{h^s \cdot \hat{\beta}(t) \cdot
  g}^{-1} \right)(U_i) \\
& = & \left( \pi_* \circ \omega_{\hat{\beta}(c(t)) \cdot g \cdot h(s,t)}^{-1} \right)(U_i) \\
& = & \left( \pi_* \circ R_{g \cdot h(s,t)*} \circ
\omega_{\hat{\beta}(c(t))}^{-1} \circ \mbox{Ad} (g \cdot h(s,t)) \right)(U_i) \\
& = & \left( \pi_* \circ \omega_{\hat{\beta}(c(t))}^{-1} \circ \mbox{Ad} (g \cdot h(s,t)
\cdot g^{-1}) \circ \Ad g \right) (U_i) \\
& = & \left( \frac{1}{1+st}\right) ^{\sigma(i)} f_i(c(t))
\end{eqnarray*}
which proves (3).
\end{Pf}



\section{Maximal degree of nilpotence implies conformal flatness}
\label{vanishing.curvature}
This section is devoted to the proof of the following proposition, the next step towards theorem \ref{flatness.thm}.
\begin{propn}
\label{prop.step.flatness}
If the group $H$ and the pseudo-Riemannian manifold $M$ satisfy the assumptions of theorem \ref{flatness.thm}, then $M$ is conformally flat.

\end{propn}

Recall that a  type $(p,q)$ pseudo-Riemannian manifold is {\it conformally flat} whenever it is locally conformally equivalent to $\Ein^{p,q}$.  If $\dim M \geq 4$, conformal flatness is equivalent to the vanishing the \emph{Weyl curvature} $W$, which is a conformally invariant $(3,1)$ tensor on $M$; if $\dim M = 3$, then vanishing of the $(3,0)$ \emph{Cotton tensor} characterizes flatness (see {\cite[p 131]{ag}}).  
%

For the canonical Cartan geometry $(M,B,\omega)$ associated to the conformal structure $(M,[\sigma])$, the Cartan curvature is defined as follows: the $2$-form
$$ \mbox{d}\omega + \frac{1}{2}[\omega, \omega]$$
on $B$ vanishes on $u \wedge v$ at $b$ whenever $u$ or $v$ is tangent to the fiber of $b$.  We define the \emph{Cartan curvature} $K$ to be the resulting function $B \rightarrow \Lambda^2 (\g / \p)^* \otimes \g$ (see \cite[5.3.22]{sharpe}). Vanishing of $K$ on $B$ is equivalent to $M$ being conformally flat (see \cite[ch 7]{sharpe} or \cite[ch IV]{kob}).  The function $K$ is $\Aut M$-invariant and $P$-equivariant; in particular, if $K(b) = 0$ for $b \in \pi^{-1}(x)$, then $K$ vanishes on the fiber of $B$ over $x$.  In this case we will also say that $K$ vanishes at $x$. 

To prove proposition \ref{prop.step.flatness}, we suppose $V \subset M$ is not flat.
Then $\partial V$ is a nonempty $H$-invariant closed subset. Under the assumptions of theorem \ref{flatness.thm}, there exists a flow $h^s$ of $H$, a point $x_0 \in \partial V$, and $b_0 \in \pi^{-1}(x_0)$ such that the holonomy of $h^s$ with respect to $b_0$ is the null translation $\tau^s$ studied  in \ref{subsection.dynamics.null}. We use the differential of the flow  $h^s$ computed above plus an idea of \cite{ccvf} to show vanishing of the Weyl and Cotton tensors along any null geodesic $\beta$ emanating from a point on $\Delta$.  Next we examine geodesic triangles in this set of vanishing curvature to show in proposition \ref{prop.vanishing.on.nbhd} that in fact the Weyl and Cotton tensors vanish in a neighborhood of $x_0$---a contradiction.  

The proof of proposition \ref{prop.vanishing.on.nbhd} below will require several preliminary results exposed in subsections \ref{sec.vanishing.lightcone} and \ref{sec.vanishing.neighborhood}.

\subsubsection{Notation}
In the two following subsections, the points $x_0$ and $b_0$ and the $1$-parameter groups $h^s$ and $\tau^s$ are as in the paragraph above.  The infinitesimal generator of $\tau^s$ is $T$ given in subsection \ref{subsection.dynamics.null}. Recall also $\lieu^-$ and the basis $U_1, \ldots, U_n$ first defined in section \ref{brief.description}.  See proposition \ref{prop.combined.reparam.framing} for the definitions of $\Delta$, $I_{\Delta}$, and $\mathcal{S} \subset \mathcal{N}(\lieu^-)$.  Recall that each curve $\beta(t) = \exp(\hat{\Delta}(v),tU)$ with $U \in \mathcal{S}$ and $v \in I_\Delta$, is defined for all $t \in \R$.

%



\subsection{Vanishing on lightcones emanating from $\Delta$}
\label{sec.vanishing.lightcone}
The aim of this subsection is the proof of:

\begin{propn}
\label{prop.cartan.vanishing}
For every $U \in \mathcal{S}$ and $v \in I_\Delta$, the Cartan curvature of $(M,B,\omega)$ vanishes on $\pi^{-1}(\beta(t))$ for all $t \in \R$, where $\beta(t) = \pi \circ \exp(\hat{\Delta}(v),tU)$.
Consequently, the Cartan curvature vanishes on the lightcone of each point of $\Delta$, in a sufficiently small neighborhood. 
\end{propn}


\begin{Pf}
Choose $v \in I_\Delta$.
We will show that when $p+q \geq 4$, the Weyl curvature vanishes on $\beta$, and the Cotton tensor vanishes when $p+q=3$.  These tensors are zero on a closed set, and $\mathcal{S} \cup - \mathcal{S}$ is dense in
 $\mathcal{N}(\lieu^-)$.  The neighborhood $V$ can be chosen to be $\pi \circ \exp_{\hat{\Delta}(v)}$, restricted to a
 neighborhood of the origin in $\lieu^-$.  Then vanishing on the entire lightcone $C(\Delta(v)) \cap V$ will follow. From the discussion
above, vanishing of the Weyl and Cotton tensors implies flatness, which implies vanishing of the Cartan curvature on
the same subset.

Let $f_i(t)$ be the framing along $\beta$ given by proposition
\ref{prop.combined.reparam.framing} (3).  We first assume $n \geq 4$ and consider
the Weyl tensor.  The conformal action of the flow $h^s$ obeys
$$ W(h^s_* f_i(t), h_*^s f_j(t), h_*^s f_k(t)) = h^s_* W(f_i(t), f_j(t), f_k(t))$$
The left hand side is
$$ \left( \frac{1}{1+st} \right)^{\sigma(i) + \sigma(j) + \sigma(k)} W(f_i(c_s(t)),f_j(c_s(t)),f_k(c_s(t)))$$

We assume $t > 0$, so that $h^s. \beta(t) \rightarrow \beta(0) = \Delta(v)$ as $s \rightarrow \infty$ by proposition \ref{prop.combined.reparam.framing} (2).  (If $t < 0$, then make $s \rightarrow - \infty$.)  Now
$$ W(f_i(0), f_j(0), f_k(0)) = \lim_{s \rightarrow \infty} (1+st)^{\sigma(i) + \sigma(j) + \sigma(k)} h^s_* W(f_i(t), f_j(t), f_k(t))$$

If $i = j = k = 1$, then the left side vanishes because $W$ is skew-symmetric in the first two entries.  Therefore, we may assume the $\sigma(i) + \sigma(j) + \sigma(k) \geq 1$.  Boundedness of the right hand side implies 
$$ h_*^s W(f_i(t), f_j(t), f_k(t)) \rightarrow 0 \qquad \mbox{as} \ s \rightarrow \infty$$
Because $h_*^s(f_1(t)) = f_1(c(t))$, the above limit means $W(f_i(t), f_j(t), f_k(t))$ cannot have a nontrivial component on $f_1(t)$.  Then
$$ W(f_i(0), f_j(0), f_k(0)) \in \mbox{span} \{ f_2(0), \ldots, f_n(0) \} $$
Varying $U$ over $\mathcal{S}$,
one sees that the Weyl curvature at $\Delta(v)$ has image in
\begin{eqnarray*}
& & \bigcap_{g \in S} \pi_* \omega_{\hat{\Delta}(v)}^{-1}(\mbox{span} \{ (\Ad g)(U_2), \ldots, (\Ad g)(U_n) \}) \\
& = & \bigcap_{g \in S} \pi_* \omega_{\hat{\Delta}(v)}^{-1} (\Ad g)(U_n)^\perp
\end{eqnarray*}

By (2) of lemma \ref{lmm.alg.facts}, the set of all $(\Ad g)(U_n)$ with $g \in S$ is a dense set of directions in the null cone $\mathcal{N}(\mathfrak{u}^-)$.  Then the intersection above is $0$,
so $W$ vanishes at $\Delta(v)$.


Now 
$$ 0 = \lim_{s \rightarrow \infty} (1+st)^{\sigma(i) + \sigma(j) + \sigma(k)} h^s_* W(f_i(t), f_j(t), f_k(t))$$
If $\sigma(i) + \sigma(j) + \sigma(k) \geq 2$, then 
$$W(f_i(t), f_j(t), f_k(t)) = 0$$  
because $h_*^s$ cannot contract any tangent vector at $\beta(t)$ strictly faster than $(1+st)^2$.  If $\sigma(i) + \sigma(j) + \sigma(k) = 1$, then we may assume $i=k=1$, and $h_*^s$ must contract the Weyl curvature strictly faster than $(1+st)$, which is possible only if
$$W(f_1(t), f_i(t), f_1(t)) \in \R f_n(t)$$
 But, in this case, for any inner product $\langle , \rangle$ in the conformal class,
$$ \langle W(f_1(t), f_i(t), f_1(t)), f_1(t) \rangle = - \langle W(f_1(t), f_i(t), f_1(t)), f_1(t) \rangle = 0$$

which implies $W(f_1(t), f_i(t), f_1(t)) = 0$, and again $W$ vanishes at
$\beta(t)$, as desired.  

When $\dim M = 3$, the argument follows the same steps and is easier.  We leave it to the reader.
%
\end{Pf}


\subsection{Vanishing on a neighborhood of $x_0$}
\label{sec.vanishing.neighborhood}
The previous subsection established vanishing of the Cartan curvature
$K$ on the union of lightcones emanating from the null geodesic segment
$\Delta$ containing $x_0$.  This union does not, however, contain a
neighborhood of $x_0$ in general.  In this subsection we will show that $\Delta$, or a particular reparametrization of it, is
complete, and that lightcones of points on $\Delta$ intersect a
neighborhood of $x_0$  in a dense subset.  Then vanishing of $K$ in a neighborhood of $x_0$ will follow.  We keep the notations of the previous section: there is a flow $h^s$ of $H$ fixing $x_0$ with holonomy the lightlike translation $\tau^s$.  Recall that $T$ denotes the infinitesimal generator of the one-parameter group $\tau^s$.



\begin{propn}
\label{prop.gtheta}
There exists $g_\theta$ in the centralizer of $T$ such that $(\Ad g_\theta)(\lieu^-)$ is transverse to $\p$ and
such that the curve $\hat{\Delta}(t) = \exp(b_0,t (\Ad g_\theta) (U_1))$ in $B$ is
defined for all time $t$.   
\end{propn}

\begin{Pf}
Recall that $x_0$ and $\tau^s$ were obtained by theorem \ref{degree.bound.precise}, which ensured that $\h^{b_0}$ centralizes $T$ (see the beginning of section \ref{section.general.degree.bound} for the notation $\h^{b_0}$).  Recall the dynamics on $\Ein^{p,q}$ of $\tau^s$ (fact \ref{fact.taus.dynamics}): for each $y$ in the null geodesic $\Lambda$, an open
dense subset of the cone $C(y)$ tends under $\tau^s$ to $y$.  Then any flow coming from
the centralizer of $T$ must leave $\Lambda$ setwise invariant; in particular, $\h^{b_0}$ preserves $\Lambda$.

\begin{lmm}
\label{lmm.2pts.stab.lambda}
Let $\lien$ be a nilpotent subalgebra of $o(p+1,q+1)$ fixing two points on $\Lambda$.
Then the nilpotence degree of $\lien$ is at most $2p$.
\end{lmm}

\begin{Pf}
The stabilizer in $\mbox{PO}(p+1,q+1)$ of a lightlike geodesic $\Lambda$ in $\Ein^{p,q}$ restricts to an action equivalent to $\PSL(2,\R)$ on $\R \BP^1$; in particular, the stabilizer acts transitively on pairs of distinct points on $\Lambda$.  Thus we may assume $\lien$ fixes $[e_0]$ and $[e_1]$.  Fixing $[e_0]$ means $\lien$
is a subalgebra of $\p \cong \co(p,q) \ltimes
\R^{p,q}$. 
Recall the inverse stereographic projection $\varphi : \R^{p,q} \rightarrow \Ein^{p,q}$ from section \ref{intuition}.
Let $u_1, \ldots, u_n$ be the standard basis of $\R^{p,q}$.  Then $\lim_{t \rightarrow \infty} \varphi(t u_1) = [e_1]$.
As in section \ref{intuition}, the set of lines in $\R^{p,q}$ which tend to
$[e_1]$ all have the form $\{ x + t u_1 \}$ for $x \in u_1^\perp$.
This set of lines is invariant by the $\lien$-action on $\R^{p,q}$, which
means that the translational components of $\lien \subset \co(p,q) \ltimes
\R^{p,q}$ are all in $u_1^\perp$, and the linear components preserve $\R u_1$,
and therefore also $u^\perp_1$.  Now by calculations similar to those in the
proof of \ref{opq.deg.bound}, we see that, if $\overline{\lien}$ is the
projection of $\lien$ on $\co(p,q)$, then
$$ \mathfrak{n}_k \subseteq \overline{\mathfrak{n}}_k + \overline{\mathfrak{n}}^k (u_1^\perp)$$
for each positive integer $k$.  But the nilpotence degree of a nilpotent subalgebra $\overline{\lien}$ of
$\co(p,q)$ is at most $2p-1$, while the order of $\overline{\lien}$ on
$u_1^\perp$ is easily seen to be at most $2p$ (compare with lemma \ref{rep.order.bounds}).
\end{Pf}

As in the proof above, the image of the restriction $\liea$
of $\h^{b_0}$ to $\Lambda$ is isomorphic to a subalgebra of $\mathfrak{sl}(2,\R)$.  Because $\liea$ is
nilpotent, $\dim \liea \leq 1$.  Let $\mathfrak{v}$ be the unipotent radical of $\mathfrak{h}^{b_0}$; it has nilpotence degree $2p+1$ since $\h^{b_0}$ does.  The restriction of $\mathfrak{v}$ to $\Lambda$ is generated by a parabolic element of $\mathfrak{sl}(2,\R)$, and it must be nontrivial by lemma \ref{lmm.2pts.stab.lambda}.  Thus $\liea$ is of parabolic type.  
Let $L \in \h^{b_0}$ have nontrivial image in $\liea$; denote this image by $\bar{L}$.  


Suppose $\bar{L}$ is parabolic type and that it fixes $[1:0] \in \R\BP^1$.  Then $L$ fixes $[e_0]$ in $\Ein^{p,q}$.  The
$1$-parameter group $e^{sL}$ preserves $\Lambda(t) = \pi(e^{tU_1})$ and reparametrizes it by $t \mapsto
\frac{t}{1+st}$.  Suppose that $\hat{\Delta}(t) = \exp(b_0,tU_1)$ is defined
on $(-\epsilon,\epsilon)$.  Take $s_\infty = -1/\epsilon$ and $s_{- \infty} =
1/\epsilon$ and apply proposition \ref{prop.holonomy} (1) to see that
$\exp(b_0,tU_1)$ is defined for all $t \in \R$.

Next suppose that $\bar{L}$ fixes $[0:1] \in
\R\BP^1$.  Then $e^{t\bar{L}}.[1:0] = [1:t]$, and, in $\Ein^{p,q}$, the orbit
is $ e^{tL}.[e_0] = \Lambda(t)$.  Then there exist $g(t) \in P$ such that
\begin{eqnarray}
\label{equation.etL}
e^{tL} = e^{tU_1} \cdot g(t)
\end{eqnarray}
 Then there exist
$(\alpha,\beta) \subset \R$, a diffeomorphism $c : (\alpha,\beta) \rightarrow
\R$, and a path $g(t) \in P$ such that 
$$ e^{tL} = e^{c(t) U_1} \cdot g(t)$$
for all $t \in (\alpha, \beta)$.  The curve $\exp(b_0,tL)$ is the orbit of $b_0$ under the lift of a conformal flow, so it is complete, and proposition \ref{prop.flow.completeness}
applies to give that $\exp(b_0,tU_1)$ is defined for all $t \in \R$.
Note that, because the two subgroups $e^{tL}$ and $e^{tU_1}$ have the same
restriction to $\Lambda$, the path $g(t)$ is in the subgroup $P_\Lambda < P$ pointwise
fixing $\Lambda$. 

Last, consider arbitrary $\bar{L}$ of parabolic type.   There exists
$\bar{g}_\theta \in \SL(2,\R)$ a rotation such that $(\Ad \bar{g}_\theta )(\bar{L})$
fixes $[0:1]$.  Let $g_\theta$ be the image of $\bar{g}_\theta$ under the standard embedding $\SL(2,\R) \rightarrow \mbox{PO}(p+1,q+1)$ given by
$$\left\{ \left( \begin{array}{ccc}
A&0 &  0   \\
0&I_{n-2} &  0   \\
  0&0 & A^{-1}  \\
  \end{array} \right)   
\ : \qquad 
\begin{array}{c}
A \in \SL(2,R) \\
  
\end{array}
\right\}
$$
with respect to which the identification $\R\BP^1 \rightarrow \Lambda$ is equivariant. Then $g_\theta$ centralizes $T$.  In $\mbox{PO}(p+1,q+1)$,
$$ g_\theta e^{tL} g_\theta^{-1} = e^{tU_1} \cdot g(t)$$
 where $g(t) \in P_\Lambda$ is as in (\ref{equation.etL}).  So
$$ e^{tL} = e^{(\Ad g_\theta)(t U_1)} \cdot h(t)$$
where $h(t) = g_\theta g(t) g_\theta^{-1}$.  The
subgroup $P_\Lambda$ is normalized by $g_\theta$, so $h(t) \in P_\Lambda$.  Proposition
\ref{prop.flow.completeness} applies to show $\exp(b_0,(\Ad g_\theta)(tU_1))$ is complete, because $\exp(b_0,tL)$ is defined for all $t$.

To prove the transversality claim, we show $(\Ad g_\theta)(\lieu^-)$ is still transverse
to $\p$, provided $g_\theta$ does not exchange $[e_0]$ and $[e_1]$ in
$\Ein^{p,q}$. Then we will take $g_\theta = 1$ when $\bar{L}$ fixes $[1:0]$, and to be the above rotation when $\bar{L}$ is parabolic but
does not fix $[1:0]$.

The subalgebra $(\Ad g_\theta)( \lieu^-)$ is transverse to $\p$ if the
orbit of $[e_0]$ in $\Ein^{p,q}$ under it is $n$-dimensional.  In the Minkowski chart
${\bf M}([e_{n+1}])$, the point $[e_0]$ is the origin, and $\Lambda$
is a null line through the origin, meeting the lightcone at infinity in one
point, $[e_1]$.  The subalgebra $\lieu^-$ acts by translations. If $g_\theta$ does not exchange $[e_0]$ and $[e_1]$, then
$g_\theta^{-1} [e_0]$ is a point on $\Lambda$ still contained in
${\bf M}([e_{n+1}])$.  The orbit 
$$(g_\theta \lieu^- g_\theta^{-1}).[e_0] = g_\theta ({\bf M}([e_{n+1}]))$$ 
which is $n$-dimensional.
\end{Pf}

\begin{propn}
\label{prop.can.conjugate}
Let $g_\theta \in \mbox{PO}(p+1,q+1)$ be given by proposition \ref{prop.gtheta}, and
$\mathcal{S}$ as in proposition \ref{prop.combined.reparam.framing}.  Let
$\hat{\Delta}(v) = \exp(b_0,v (\Ad g_\theta) (U_1))$ and $\Delta = \pi \circ
\hat{\Delta}$.  Let $\mathcal{S}' =
(\Ad g_\theta)(\mathcal{S})$.  Then 
\begin{enumerate}
\item{The flow $h^s$ fixes $\Delta(v)$ pointwise and has holonomy at $\Delta(v)$ with respect to $\hat{\Delta}$ equal $\tau^s$.}
\item{For each $U \in \mathcal{S}'$ and $v \in \R$, the curve
    $\hat{\beta}(t) = \exp(\hat{\Delta}(v),tU)$ is complete and projects to a null geodesic.}
\item{For $\Delta$ as in (1) and $\hat{\beta}$ as in (2), the Cartan curvature
    vanishes on the fiber of $\hat{\beta}(t)$ for all $t \in \R$.}
\end{enumerate}
\end{propn}

\begin{Pf}
\begin{enumerate}
\item{Since $g_\theta$ centralizes the null translation $T$, one can follow the same proof as for proposition \ref{prop.combined.reparam.framing} (1).}
\item{Recall the embedding $\PSL(2,\R) \rightarrow \mbox{PO}(p+1,q+1)$ with image the stabilizer of $\Lambda$.  The element $g_\theta$ corresponds to rotation by some angle $\theta$ on $\Lambda = \BP(\mbox{span}\{ e_0,e_1 \})$.  Compute that for 
$U_n = E_n^0 - E_{n+1}^1$ as in lemma \ref{lmm.alg.facts},
$$(\Ad g_\theta)(U_n) = U_n$$
Next, note that the subgroup $S$ of lemma \ref{lmm.alg.facts} is contained in
$P_\Lambda$, the pointwise stabilizer of $\Lambda$.  Then
$$ g_\theta S g_\theta^{-1} < P_\Lambda < P$$
Now any element of $\mathcal{S}'$ is of the form 
$$\lambda (\Ad g_\theta \circ \Ad w)(U_n)$$
for some $\lambda \in {\bf R}$ and $w \in S$, and can be written
$$ \lambda (\mbox{Ad}(g_\theta w g_\theta^{-1}))(U_n)= U$$
Since $g_\theta w g_\theta^{-1} \in P$, the element $U$ projects to a
null vector in $\g/\p \cong \R^{p,q}$.  Then $\pi \circ
\exp(\hat{\Delta}(v),tU)$ is a null geodesic.

Fix $U = (\Ad g_\theta w) (U_n) \in \mathcal{S}'$; it suffices to prove (2) for such $U$, since the geodesic generated by $\lambda U$ is complete if and only if the geodesic generated by $U$ is.  Recall the matrices $h(s,t)$, representing
the holonomy of $h^s$ along null geodesics based at $\pi(\exp(b_0,vU_1))$ with initial direction in
$\mathcal{S}$ or $- \mathcal{S}$.  Straightforward computation shows that $g_\theta$ commutes
with $h(s,t)$.  Recall also that $e^{sT}$ commutes with $S$.  Then we compute in $\mbox{PO}(p+1,q+1)$,
\begin{eqnarray*}
e^{sT} e^{tU} 
 & = & e^{c(t)U} (g_\theta w g_\theta^{-1}) h(s,t) (g_\theta w
g_\theta^{-1})^{-1} 
\end{eqnarray*}
where $c(t) = \frac{t}{1+st}$.

Now proposition \ref{prop.holonomy} with part (1) implies that
$h^s$ reparametrizes $\beta(t)$ by $c(t)$ and has holonomy 
$$g_\theta w h(s,t) w^{-1} g_\theta^{-1} = (g_\theta w g_\theta^{-1}) h(s,t) (g_\theta w
g_\theta^{-1})^{-1} \in P$$
along it with respect to $\hat{\beta}$.  Part (1) of proposition
\ref{prop.holonomy} gives the desired completeness.
}

\item{ As above, we may assume $U = \mbox{Ad}(g_\theta w)(U_n)$.  Define a framing along $\beta(t)$ as in proposition
\ref{prop.combined.reparam.framing} by 
$$f_i(\beta(t)) = (\pi_* \circ \omega_{\hat{\beta}(t)}^{-1} \circ (\Ad g_\theta w)) (U_i)$$
Recall that $(\Ad g_\theta)(\lieu^-)$ is
transverse to $\p$ by proposition \ref{prop.gtheta}, so $\mbox{Ad}(g_\theta w)(\lieu^-)$ is, as well.  Now the derivative of
$h(s,t)$ along $\beta(t)$ in this framing is computed as in
the proof of \ref{prop.combined.reparam.framing} by the adjoint action of the holonomy $g_\theta
w h(s,t) w^{-1} g_\theta^{-1}$ on the $(\Ad g_\theta w)(U_i)$, modulo $\p$.  The derivative has the same diagonal form
as in proposition \ref{prop.combined.reparam.framing}.  In fact, all the
conclusions of \ref{prop.combined.reparam.framing}, and thus also the arguments of
proposition \ref{prop.cartan.vanishing}, hold when $\hat{\Delta}(v) =
\exp(b_0,(\Ad g_\theta)(vU_1))$, and $\mathcal{S}$ is replaced
by $\mathcal{S}'$, so we conclude that the Cartan curvature vanishes
along the desired geodesics.}
\end{enumerate}
\end{Pf}

\begin{defn}
Let $(M,B,\omega)$ be a Cartan geometry modeled on $G/P$.  Let $\gamma$ be a piecewise smooth
curve in $B$.  The \emph{development} $\mathcal{D} \gamma$ is the piecewise smooth curve in $G$
satisfying $\mathcal{D}\gamma (0) = e$ and $(\mathcal{D} \gamma)'(t) =
\omega(\gamma'(t))$ for all but finitely many $t$.
\end{defn}

Note that for any piecewise smooth curve $\gamma$ in $B$, the development $\mathcal{D} \gamma$ is
defined on the whole domain of $\gamma$, because it is given by a linear first-order ODE on $G$ with bounded coefficients.   

%

\begin{propn}
\label{prop.vanishing.on.nbhd}
The Cartan curvature $K$ vanishes on an open set of the form $\pi^{-1}(V)$ for
$V$ a neighborhood of $x_0$ in $M$.
\end{propn}

\begin{Pf}
Recall the basis $U_1, \ldots, U_n$ for $\lieu^-$ first introduced in section \ref{brief.description}. 
Let $\hat{\Delta}(v) = \exp(b_0,(\Ad g_\theta)(vU_1))$, where $g_\theta$ is
given by proposition \ref{prop.gtheta}, so $\hat{\Delta}$ is complete.  Recall
that the curves $\exp(\hat{\Delta}(v),tU)$, where $U \in (\Ad
g_\theta)(\mathcal{S}) = \mathcal{S}'$ are complete, as well, from proposition
\ref{prop.can.conjugate}.

To show that the Cartan curvature vanishes on a neighborhood above $x_0$, it
suffices to show that $K=0$ on $\exp(b_0,V)$, for $V$ a neighborhood of $0$ in
$(\Ad g_\theta)(\lieu^-)$, because $\pi_{*b_0}$ maps $\omega_{b_0}^{-1}(V)$ onto
a neighborhood of $0$ in $T_{x_0}M$ by proposition \ref{prop.gtheta}.  

First suppose $g_\theta = 1$.  Recall from the proof of lemma
\ref{lmm.alg.facts} (2) that $\mathcal{S}$ consists of all $U \in \mathcal{N}(\lieu^-)$ with $\langle
U,U_1 \rangle > 0$.  Let 
$$Y = a U_1 + X + c U_n \in \lieu^-$$
with $X \in \mbox{span} \{ U_2, \ldots, U_{n-1} \}$, and assume that $c \neq
0$.  Let $b = \langle X, X \rangle$.  Define, for each $0 \leq r \leq 1$, a piecewise smooth
curve $\alpha_r$ in $B$ by concatenating 
$$\hat{\Delta}(t(b^2/c + 2a)), \ 0 \leq t \leq r/2$$
 and 
$$\exp(\hat{\Delta}(rb^2/2c + ra), (2t-r)(cU_n + X - b^2/2c \ 
U_1)), \ r/2 \leq t \leq r$$
Note that $cU_n + X - b^2/2c\ U_1 \in \pm \mathcal{S}$ because $c \neq 0$.  Define $\beta(r) = \alpha_r(r)$.  

\begin{figure}[ht]
{ \input{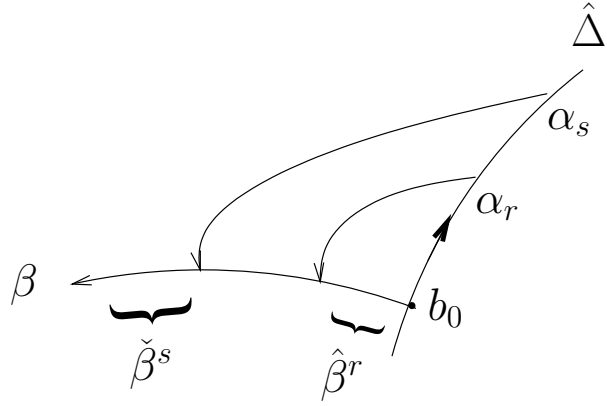}}
\caption{components of the homotopy between $\beta$ and $\alpha_1$}
\end{figure}

Because $\lieu^-$ is an abelian subalgebra of $\g$, the development
$$\mathcal{D}\alpha_r(r) = e^{r(\frac{b^2}{2c}+a)U_1} \cdot e^{r(cU_n + X - \frac{b^2}{2c} U_1)} = e^{rY}$$
  Denote by
$\hat{\beta}^s$ the restriction of $\beta$ to $[0,s]$ and by $\check{\beta}^s$ the restriction of $\beta$ to $[s,1]$.   
The curve $\beta$ is homotopic to $\alpha_1$ through the family of
concatenations $\alpha_s * \check{\beta}^s$; similarly, $\hat{\beta}^r$ is
homotopic to $\alpha_r$ for all $0 \leq r < 1$. 
Because the curvature $K$ vanishes on the images of these homotopies, the developments of the homotopic curves have the same endpoints (see \cite[3.7.7 and 3.7.8]{sharpe}):
$$\mathcal{D}\hat{\beta}^r(r) = \mathcal{D}\alpha_r(r) = e^{rY} \ \forall \ r \in [0,1]$$
But 
$$\mathcal{D}\beta(r) = \mathcal{D}\hat{\beta}^r(r)= e^{rY}$$
 which
means that
$$ \beta(r) = \exp(b_0,rY) $$

Then $K(\exp(b_0,Y)) = 0$.  Varying $Y$ over all sufficiently small $aU_1
+ X + cU_n$ with $c \neq 0$ and passing to the closure gives vanishing of $K$
on $\exp(b_0,V)$, for $V$ a neighborhood of $0$ in $\lieu^-$.

If $g_\theta \neq 1$, then consider 
$$Y  = (\Ad g_\theta)(a U_1 + X + c U_n) \in (\Ad g_\theta)(\lieu^-)$$
again with $c \neq 0$.  Define $\alpha_r$ by concatenating a portion of
$\hat{\Delta}$ as above with
$$ \exp(\hat{\Delta}(rb^2/2c + ra),(2t-r) (\Ad g_\theta)(cU_n + X - b^2/2c\  U_1)), \
r/2 \leq t \leq r$$
note that $(\Ad g_\theta)(cU_n + X - b^2/2c U_1) \in \pm \mathcal{S}'$.
Therefore, by proposition \ref{prop.can.conjugate} (2) and (3), the curves
$\alpha_r$ are defined on $[0,r]$ for all $r$ and $K$ vanishes on them.  Then
the same argument as above applies to give vanishing of curvature at
$\exp(b_0,Y)$.  The set of possible $Y$ are dense in a neighborhood $V$ of $0$
in $(\Ad g_\theta)(\lieu^-)$, so we obtain the desired vanishing on $\exp(b_0,V)$.  
\end{Pf}



\section{End of the proof of theorem \ref{flatness.thm}: global structure of $M$}
\label{global}
This section is again under the assumptions of theorem \ref{flatness.thm}: $H$ is a connected nilpotent Lie group acting conformally on a type $(p,q)$ pseudo-Riemanniann manifold $(M,\sigma)$, and $d(H)=2p+1$ with $p+q \geq 3$.  From proposition \ref{prop.step.flatness}, the manifold $(M, \sigma)$ is locally conformally equivalent to $\Ein^{p,q}$, or in other words is endowed with a $(\widetilde{G},{\widetilde \Ein}^{p,q})$-structure, where $\widetilde{G} = \Conf \widetilde{\Ein}^{p,q}$, a covering group of $\mbox{PO}(p+1,q+1)$.  The {\it developing map} of the structure is a conformal immersion of the universal cover of $M$
\[   \delta : {\widetilde M} \to {\widetilde \Ein}^{p,q} \]
(see \cite{thurston}, \cite{goldman} for an introduction to $(G,X)$-structures and the construction of the developing map).  
If $({\widetilde M},{\widetilde B},\tilde \omega)$ is the canonical Cartan geometry associated to the lifted conformal structure $({\widetilde M}, \tilde{\sigma})$, then $\delta$ lifts to an immersion of bundles $\hat \delta : {\widetilde B} \to \widetilde{G}$, where $\widetilde{G}$ is seen as a principal bundle over $\widetilde{\Ein}^{p,q}$.  This immersion satisfies ${\hat \delta}^*\omega_{\widetilde{G}}=\tilde \omega$, where $\omega_{\widetilde{G}}$ denotes the Maurer-Cartan form on $\widetilde{G}$.  

Let $\rho$ be the {\it holonomy morphism} $\Conf {\widetilde M} \to {\widetilde G}$, related to the developing map by the equivariance property
\[   \rho(\phi) \circ \delta = \delta \circ \phi \ \ \ \forall\ \phi \in \Conf {\widetilde M}.\]


The $H$-action lifts to a faithful action of a connected covering group, which we will also denote $H$, on ${\widetilde M}$. 
The group $\rho(H)=\check H$ is a connected nilpotent subgroup of $\widetilde{G}$, with Lie algebra $\check \h$ isomorphic to $\h$.  In particular $d(\check \h)=2p+1$, and by proposition \ref{opq.deg.bound}, we may assume that $\check H$ contains the lifted $1$-parameter group $\tau^s$ (see section \ref{subsection.dynamics.null}).   

The fundamental group $\pi_1(M)$ is isomorphic to a discrete subgroup $\Gamma < \Conf {\widetilde M}$ acting freely and properly on ${\widetilde M}$.  
 Because $\Gamma$ centralizes $\h$, the image $\rho(\Gamma) = \Phi$ centralizes $\check \h$.

To complete the proof of theorem \ref{flatness.thm}, we must prove that {\it the developing map $\delta$ is a conformal diffeomorphism between ${\widetilde{M}}$ and ${\widetilde \Ein}^{p,q}$}.  Subsections \ref{more.geometry}, \ref{sec.centralizer}, and \ref{geometrical.properties} contain preliminary  geometric and algebraic results to this end.  The end of the proof of theorem \ref{flatness.thm} will be given in \ref{sec.proof.lorentzian} for the Lorentzian case, and in \ref{sec.proof.p}  for the other types.

\subsection{More on geometry and dynamics on $\tein$}
\label{more.geometry}
When $p\geq 2$, ${\widetilde \Ein}^{p,q}$ is a double cover of $\Ein^{p,q}$.  It is $\widehat{\mathcal{N}}^{p+1,q+1}/\R_{>0}^*$. 
The conformal group $\widetilde{G}$ is $\mbox{O}(p+1,q+1)$, and the stabilizer of $[e_0]$ is an index-two subgroup of $P$. 
A lightcone $C(x)$ in ${\widetilde \Ein}^{p,q}$ has two singular points, and its complement has two connected components, each one conformally equivalent to $\R^{p,q}$.

The Lorentz case $p=1$ is more subtle since ${\widetilde \Ein}^{1,n-1}$ is no longer compact.  It is conformally equivalent to $(\R \times {\bf S}^{n-1},-dt^2 \oplus g_{{\bf S}^{n-1}})$.  Details about this space are in \cite[ch 4.2]{charlesthese} and \cite{primer}. The group $\widetilde{G} = \Conf {\widetilde \Ein}^{1,n-1}$ is a twofold quotient of the universal covering group of ${\mbox O}(2,n)$, with center $Z \cong {\bf Z}$.  The space $\Ein^{1,n-1}$ is the quotient of ${\widetilde \Ein}^{1,n-1}$ by the $Z$-action. 

The lightlike geodesics and lightcones in ${\widetilde \Ein}^{1,n-1}$ are no longer compact.  Any lightlike geodesic can be parametrized $\gamma(t) = (t,c(t))$, where $c(t)$ is a unit-speed geodesic of ${\bf S}^{n-1}$.  Lightlike geodesics are preserved by $Z$, which acts on them by translations; the quotient is a lightlike geodesic of $\Ein^{1,n-1}$.  Any lightcone $C(x) \subset {\widetilde \Ein}^{1,n-1}$ has infinitely-many singular points, which coincide with  the $Z$-orbit of $x$.  The complement of $C(x)$ in ${\widetilde \Ein}^{1,n-1}$ has a countable infinity of connected components, each one conformally diffeomorphic to $\R^{1,n-1}$.  The center $Z$ freely and transitively permutes these Minkowski components.

Recall that $\tau^s$ is the flow on $\widetilde{\Ein}^{1,n-1}$ generated by the null translation $T$.  The null geodesic $\Lambda = \BP(\mbox{span}\{e_0,e_1 \})$ is the fixed set of $\tau^s$ on $\Ein^{1,n-1}$.  Let $\tlambda$ be the inverse image of $\Lambda$ in $\widetilde{\Ein}^{1,n-1}$; it is noncompact and connected, and equals the fixed set of $\tau^s$ on $\widetilde{\Ein}^{1,n-1}$.  Given $\tx \in {\widetilde \Ein}^{1,n-1} \setminus \tlambda$, there are two distinct points $\tx^+$ and $\tx^-$ on $\tlambda$ such that 
$$\lim_{s \to \infty} \tau^s.\tx=\tx^+ \qquad \mbox{and} \qquad \lim_{s \to -\infty} \tau^s.\tx=\tx^-$$
Details about this material can be found in \cite[p 67]{charlesthese}.

\subsection{About the centralizer of $\check \h$}
\label{sec.centralizer}
For arbitrary $(p,q)$, let $\tlambda$ be the inverse image in $\tein$ of $\Lambda$; it is connected and fixed by $\tau^s$.  Our first task  is to find  an algebraic restriction on $\Phi$, using that it commutes with $\check{\h}$.
\begin{propn}
\label{finite.action}
The centralizer $C(\check \h)$ of $\check \h$ in $\widetilde{G}$ leaves $\tlambda$ invariant.  The $C(\check \h)$-action on $\tlambda$ factors through a homomorphism to $\Z_2$ if $p \geq 2$, and it factors through a homomorphism to $Z$ when $p=1$.
\end{propn}

\begin{Pf}
Denote $\mathfrak{c}(\check \h)$ the Lie algebra of the centralizer of $\check \h$.  Observe first that both $\mathfrak{c}(\check \h)$ and $\check \h$ centralize $\tau^s$.  Recall that $\tilde{\Lambda}$ is a local attracting set for $\tau^s$ (see the proof of proposition \ref{prop.can.conjugate}), so $\mathfrak{c}(\tau^s)$ leaves $\tilde{\Lambda}$ invariant.  The centralizer $\mathfrak{c}(\tau^s)$ in $\oo(p+1,q+1)$ of $\{ \tau^s \}$ thus consists of matrices of the form
\begin{eqnarray}
\label{eqn.ctaus.matrix}
& \qquad &
\left( \begin{array}{ccccc}
a &b &  -x^t.J_{p-1,q-1}&s  &0  \\
c&-a &  -y^t.J_{p-1,q-1} & 0  &-s  \\
  & & M &y& x \\
  & & &   a  & -b\\
  & & & -c & -a
\end{array} \right) 
\qquad 
\begin{array}{c}
a,b,c,s \in \R \\
  x,y\in \R^{p-1,q-1} \\
  M \in \oo(p-1,q-1) 
\end{array}
\end{eqnarray}

As in the proof of proposition \ref{prop.can.conjugate}, the projection of $\check \h$ on ${\mathfrak{sl}(2,\R)}$ is a $1$-dimensional subalgebra of parabolic type.  The projection of $\mathfrak{c}(\check \h)$ lies in the same subalgebra; we wish to show it is zero.




We may assume, by conjugating $\check{\h}$ in the centralizer of $\tau^s$ if necessary, that $\mathfrak{c}(\check \h)$ and $\check \h$ are subalgebras of 
$$ \mathfrak{q} = \{ X \in \mathfrak{c}(\tau^s) \ : \ a, c = 0 \}$$
This algebra is isomorphic to $(\R \oplus \oo(p-1,q-1)) \ltimes \mathfrak{heis}(2n-3)$.  Elements of $\q$ are denoted $u = (b,M,x,y,s)$,
with $b \in \R$, $x,y\in \R^{p-1,q-1} $, and $M \in \oo(p-1,q-1)$.  Denote $b=\pi_1(u)$, $M=\pi_2(u)$, and $(x,y)=\pi_3(u)$.  Note that if $\pi_i(u) = 0$ for $i = 1,2,3$, then $u$ is in the center of $\mathfrak{q}$.  If $u_1 = (b_1,M_1,x_1,y_1,s_1) $ and $u_2 = (b_2,M_2,x_2,y_2,s_2) $ are in $\q$, then an easy computation yields
\begin{itemize}
\item{$\pi_1([u_1,u_2])=0$}
\item{$\pi_2([u_1,u_2])=[M_1,M_2]$}
\item{$\pi_3([u_1,u_2])=(b_1y_2-b_2y_1-M_1.x_2+M_2.x_1,-M_1.y_2+M_2.y_1)$}
\end{itemize}
Now, if $u_0=(b_0,M_0,x_0,y_0,s_0)$ is in $\mathfrak{c}(\check \h)$, each $u=(b,M,x,y,s) \in \check \h$ must satisfy the relations:
\begin{itemize}
\item{$[M_0,M]=0$}
\item{$M_0.y=M.y_0$}
\item{$b_0y-M_0.x=by_0-M.x_0$}
\end{itemize}

We claim that $b_0  =0$.  If not, then from the last relation above, whenever $u_1$ and $u_2$ are in $\check \h$, then
$$ y_1=\frac{b_1}{b_0}y_0 - \frac{1}{b_0}M_1.x_0 + \frac{1}{b_0}M_0.x_1 $$
and 
$$ y_2=\frac{b_2}{b_0}y_0 - \frac{1}{b_0}M_2.x_0 + \frac{1}{b_0}M_0.x_2 $$
This implies
\begin{eqnarray*}
\pi_3([u_1,u_2])& = & 
(\frac{b_1}{b_0}(-M_2.x_0 + M_0.x_2)+\frac{b_2}{b_0}(-M_0.x_1+M_1.x_0)\\
& -& M_1.x_2+M_2.x_1, -M_1.y_2+M_2.y_1)
\end{eqnarray*}

As in section \ref{section.proof.deg.bound}, we may assume without affecting the nilpotence degree and without changing $C(\check{\h})$ that $\check \h$ is Zariski closed.  Write ${\check \h} \cong {\check \lier} \ltimes {\check \lieu} $, where ${\check \lieu}$ is an algebra of nilpotents, and $d(\check{\h}) = d(\check{\lieu})$. 

Let $\check{\mathfrak{m}} = \pi_2(\check \lieu)$.  It is a nilpotent Lie subalgebra of $\oo(p-1,q-1)$, and also an algebra of nilpotents, since ${\check \lieu}$ is so.  
Using the equation above, we get by induction that $\pi_3(\check \lieu_k) \subset \check{\mathfrak{m}}^k(\R^{p-1,q-1}) \times \check{\mathfrak{m}}^k(\R^{p-1,q-1})$.  Moreover, $\pi_2(\check \lieu_k)=\check{\mathfrak{m}}_k$, and $\pi_1(\check \lieu_k)=0$ as soon as $k \geq 1$.  
By proposition \ref{opq.deg.bound}, $d(\check{\mathfrak{m}}) \leq 2p-3$, and $o(\check{\mathfrak{m}}) \leq 2p-1$ by lemma \ref{rep.order.bounds}. As a consequence, $\pi_1(\check \lieu_{2p-1})=\pi_2(\check \lieu_{2p-1})=\pi_3(\check \lieu_{2p-1})= 0$, which implies that $\check \lieu_{2p-1}$ is in the center of $\check \lieu$, and finally $d(\check \lieu) \leq 2p$.  Since $d(\check \h)=d(\check \lieu)$, we get  a contradiction.  Therefore, $\mathfrak{c}(\check \h)$ is actually the subalgebra of $\mathfrak{c}(\tau^s)$ with $b=0$ in (\ref{eqn.ctaus.matrix}).






 The projection of $C(\check{h})$ to $\mbox{PO}(p+1,q+1)$ commutes with the flow of parabolic type on $\Lambda$ generated by $\check{\h}$.  The centralizer of such a $1$-parameter subgroup in $\mbox{PSL}(2,\R)$ is itself; in particular, any finite subgroup of the centralizer is trivial.  From the above calculation, on the other hand, the identity component of $C(\check \h)$ acts trivially on $\tlambda$.  It follows that the image of $C(\check{\h})$ in $\mbox{PO}(p+1,q+1)$ acts trivially on $\Lambda$, so the restriction of $C(\check{\h})$ to $\tlambda$ factors through a homomorphism to the kernel of $\widetilde{G} \rightarrow \mbox{PO}(p+1,q+1)$.  From the previous section, this kernel is isomorphic to $\Z_2$ when $p \geq 2$ and equals $Z$ when $p=1$.  This concludes the proof of the proposition. \end{Pf}

\subsection{Geometric properties of the developing map}
\label{geometrical.properties}

Recall that $\{ h^s \} < H$ satisfies $\rho(h^s)=\tau^s$.  
We now adopt the notations and results of subsection \ref{subsection.dynamics.on.M}. The lightlike geodesic $\Delta(t) = \pi \circ \exp(b_0,tU_1)$ is pointwise fixed by $h^s$ and is locally an attracting set for it (see proposition \ref{prop.combined.reparam.framing} (2)).  Choose $\tilde{x}_0 \in \widetilde{M}$ over $x_0$ and lift $h^s$ to $\widetilde{M}$.  Let $\tlambda \subset \tein$ be as above: it is a closed subset, pointwise fixed by $\tau^s$, and it is the attracting set for $\tau^s$.  Because of these dynamics, $\pi_M^{-1}(\Delta) \subset \delta^{-1}(\tlambda)$, which is a closed, $\Gamma$-invariant, $1$-dimensional, immersed submanifold.  Let $\tdelta$ be the component of $\delta^{-1}(\tlambda)$ containing $\tilde{x_0}$.  Denote $\Gamma_0 < \Gamma$ the subgroup leaving $\tdelta$ invariant.  


\begin{propn}
\label{proposition.more.descriptive.than.charles.labels}
The image $\pi_M(\tdelta) \subset M$ is closed.  Therefore, $\Gamma_0$ acts cocompactly on $\tdelta$.
\end{propn}

\begin{Pf}
We show $\pi_M(\tdelta)$ is closed in $\pi_M(\delta^{-1}(\tlambda))$, and therefore in $M$.  Suppose $\pi_M(\tilde{x}_n) \rightarrow \pi_M(\tilde{y})$ with $\tilde{x}_n \in \tdelta$ and $\tilde{y} \in \delta^{-1}(\tlambda)$.  Let $U$ be a neighborhood of $\tilde{y}$ that maps diffeomorphically to its images under $\pi_M$ and under $\delta$.  There exist $\gamma_n \in \Gamma$ such that $\gamma_n .\tilde{x}_n \to \tilde{y}$ in $U$.  Then $\delta(\gamma_n.\tilde{x}_n) \to \delta(\tilde{y})$ in $\tlambda$, and we may assume $U$ is small enough that $\delta(U) \cap \tlambda$ is an open segment.  Then $\gamma_n.\tilde{x}_n$ and $\tilde{y}$ are in a common segment of $\delta^{-1}(\tlambda) \cap U$.  Then for some $\gamma_n = \gamma$, the translate $\gamma.\tilde{y} \in \tdelta$. \end{Pf}

As a consequence, we obtain the following result.
\begin{propn}
\label{onto}
The map $\delta$ is a covering map from $\tilde \Delta$ onto $\tilde \Lambda$. When $M$ is Lorentzian, $\delta$ is a diffeomorphism between $\tilde \Delta$ and $\tilde \Lambda$.
\end{propn}

\begin{Pf}
First note that $\tdelta$ is open in $\delta^{-1}(\tlambda)$.  For if $\delta^{-1}(\tlambda)$ were recurrent, then $\tlambda$ would be, as well; but $\tlambda$ is a closed, embedded submanifold of $\tein$.  Therefore the image $\delta(\tdelta)$ is a connected open subset of $\tlambda$.  By equivarience of $\delta$ and the previous proposition, $\rho(\Gamma_0)$ preserves $\delta(\tdelta)$ and acts cocompactly on it.  But $\rho(\Gamma_0)$ centralizes $\check \h$, so its action on $\tlambda$ factors through either a finite group, or the extension of a finite group by $\Z$.  In both cases, $\delta(\tdelta)$ must equal $\tilde \Lambda$.

When $M$ has Lorentz type, then $\delta: \tdelta \to \tlambda$ must be a diffeomorphism, because all lightlike geodesics in $\widetilde{\Ein}^{1,n-1}$ are embedded copies of $\R$, as described in section \ref{more.geometry}; in particular, they have no self-intersection.

Now assume $p \geq 2$.  On one hand, $\Gamma_0$ acts cocompactly on $\tdelta$; on the other hand, the action of $\rho(\Gamma_0)$ on $\tlambda$ factors through a finite group.  Let $\Gamma_0^{\prime} \lhd \Gamma_0$ be such that $\rho(\Gamma_0^{\prime})$ is the kernel in $\rho(\Gamma_0)$ of restriction to $\tlambda$.  Then the restriction of $\delta$ factors 
$$
\begin{array}{cccc}
\left.\delta\right|_{\tdelta} : & \tdelta & \to & \tlambda \\
                               &   & \searrow    &  \uparrow \\
                               &          &     & \tdelta / \Gamma'_0  
\end{array}
$$
Because $\Gamma'_0$ acts freely and properly on $\tdelta$, the quotient map $\tdelta \to \tdelta / \Gamma'_0$ is a covering, and because $\Gamma'_0$ has finite index in $\Gamma_0$, its action on $\tdelta$ is cocompact.
The map $\tdelta / \Gamma'_0 \to \tlambda$ is surjective because $\left. \delta \right|_{\tdelta}$ is; it is a local diffeomorphism because $\left. \delta \right|_{\tdelta}$ and $\tdelta \to \tdelta / \Gamma'_0$ are.  By compactness of $\tdelta / \Gamma'_0$, it follows that $\tdelta / \Gamma'_0 \to \tlambda$, hence  $\tdelta \to \tlambda$, is a covering, as desired.  
\end{Pf}


\subsection{Proof of theorem \ref{flatness.thm}: the case of Lorentz manifolds}
\label{sec.proof.lorentzian}
Suppose now that $p=1$.  Let
$$ \Omega = \{ \tz \in {\widetilde M} \setminus \tdelta \ | \ \lim_{s \to \infty} h^s.\tz \ \mbox{exists and is in } \tdelta \}$$

\begin{propn}
\label{ouverts}
The set $\Omega$  is nonempty and open.  It is mapped diffeomorphically by $\delta$ onto $\widetilde{\Ein}^{1,n-1} \setminus \tlambda$.
\end{propn}

\begin{Pf}
Let us first check that $\Omega$ is nonempty. 
Recall $\tdelta$ is pointwise fixed by $h^s$.  Let $\tz_{\infty} \in \tdelta$, and choose $\tb_{\infty} \in \widetilde{B}$ above $\tz_{\infty}$ such that the holonomy of $h^s$ with respect to $\tb_{\infty}$ is $\tau^s$. Let ${\mathcal S}$ be as in proposition \ref{prop.combined.reparam.framing} and $U \in {\mathcal S}$. Consider the geodesic $\beta(t)=\pi \circ \exp(\tb_{\infty},tU)$.
It is complete by proposition \ref{prop.combined.reparam.framing} (2); further, for $t > 0$, 
$$\lim_{s \to \infty} h^s.\beta(t)=\beta(0)=\tz_{\infty}  $$ 
Then $\beta(t) \in \Omega$ for $t > 0$.

To prove that $\Omega$ is open, choose $\tz_0 \in \Omega$.  There exists $\tz_{\infty} \in \tdelta$ such that $\lim_{s \to \infty} h^s.\tz_0=\tz_{\infty}$. Since the orbits of $\tau^s$ are lightlike geodesics in $\widetilde{\Ein}^{1,n-1} $, the same is true for the orbits of $h^s$ on ${\widetilde M}$. 
Then $\tilde{z}_0$ lies on some lightlike geodesic emanating from $\tilde{z}_\infty$.  Any such geodesic not fixed by $h^s$ has the form $\pi \circ \exp(\tilde{b}_\infty,tU)$ with $U \in \mathcal{S}$.
Then for $s_0 >0$ there exist $\epsilon > 0$ and a diffeomorphism $c : (s_0, \infty) \to (0, \epsilon)$ such that, for every $ s \in (s_0,\infty)$,
$$\pi \circ \exp(\tb_{\infty},c(s)U)=h^s.\tz_0.$$
 There are a neighborhood $I$ of $\tz_{\infty}$ in $\tdelta$, a segment $\widetilde{I} \subset \widetilde{B}$ lying over $I$, and an open neighborhood  ${\mathcal U}$ of $0$ in ${\mathfrak u}^-$ such that the map 
\begin{eqnarray*}
\mu & : & I \times ({\mathcal U} \cap {\mathcal S}) \to {\widetilde M} \\
  &  & (\tz,u) \mapsto \pi \circ \exp(\tb,u)
\end{eqnarray*}
where $\tilde{b} \in \widetilde{I}$ lies over $\tilde{z}$, is defined and is a submersion.  Choosing  $s_0$ big enough, $c(s_0)U \in {\mathcal U} \cap {\mathcal S}$, so that $V=\mu(I \times ({\mathcal U} \cap {\mathcal S}))$ is an open subset containing $h^{s_0}.\tz_0$. It follows immediately from proposition \ref{prop.combined.reparam.framing} that $V \subset \Omega$.  Since $\Omega$ is $h^s$-invariant,  $h^{-s_0}(V) \subset \Omega$. It is an open subset containing $\tz_0$, which shows  that $\Omega$ is open.

We now  prove that $\delta$ is an injection in restriction to $\Omega$.  Assume that $\tz$ and $\tz^{\prime}$ are two points of $\Omega$ satisfying $\delta(\tz)=\delta(\tz^{\prime})$.   
Let $\tz_\infty = \lim_{s \to \infty} h^s.\tz$ and $\tz_{\infty}^{\prime} = \lim_{s \to \infty}h^s.\tz^{\prime}$.  Then 
$$\delta(\tz_{\infty})=\lim_{s \to \infty}\tau^s.\delta(\tz)=\lim_{s \to \infty}\tau^s.\delta(\tz^{\prime})=\delta(\tz_{\infty}^{\prime})$$
Because $\delta$ is injective on $\tdelta$ by proposition \ref{onto}, $\tz_{\infty}=\tz_{\infty}^{\prime}$.  Choose $U$ an open neighborhood of $\tz_{\infty}$ which is mapped diffeomorphically by $\delta$ on an open neighborhood $V$ of $z_{\infty}=\delta(\tz_{\infty})$.   There exists $s_0 \geq 0$ such that for all $t \geq s_0$, $h^t.\tz \in U$ and $h^t.\tz^{\prime} \in U$.  Moreover, $\delta(h^t.\tz)=\delta(h^t.\tz^{\prime})=\tau^t.\delta(\tz)$.  Since $\delta$ is an injection in restriction to $U$, the images $h^t.\tz = h^t.\tz^{\prime}$, 
so $\tz=\tz^{\prime}$, as desired.

It remains to show that $\delta(\Omega)=\widetilde{\Ein}^{1,n-1} \setminus \tlambda$.  The inclusion  $\delta(\Omega) \subset \widetilde{\Ein}^{1,n-1} \setminus \tlambda$ follows easily from the definition of $\Omega$.  Just note that any $\tz \in \delta^{-1}(\tlambda)$ is fixed by $h^s$, so $\Omega$ cannot meet $\delta^{-1}(\tlambda)$.  Now, pick $z \in \widetilde{\Ein}^{1,n-1}$.  There exists $z_{\infty} \in \tlambda$ such that $\lim_{s \to \infty}\tau^s.z=z_{\infty}.$ 
By proposition \ref{onto}, there is a unique $\tz_{\infty} \in \tdelta$ such that $\delta(\tz_{\infty})=z_{\infty}$.  Also, there is a neighborhood $U$ of $\tz_{\infty}$ mapped diffeomorphically by $\delta$ on some neighborhood $V$ of $z_{\infty}$.  There exists $s_0$ such that for $s \geq s_0$, $\tau^s.z \in V$.  Let $\tz \in U$ be such that $\delta(\tz)=\tau^{s_0}.z$.  Then for all $s \geq s_0$, we have $h^s.\tz \in U$ and $\lim_{s \to \infty} h^s.\tz=\tz_{\infty}$.  Thus, $\tz \in \Omega$.  Moreover,  $\delta(h^{-s_0}.\tz)=z$ and since $\Omega$ is $h^s$-invariant, $z \in \delta(\Omega)$, as desired.  \end{Pf}

\begin{rmk}
\label{adherence}
Notice that when we proved that $\Omega$ is nonempty, we showed that $\tdelta$ is in the closure of $\Omega$.
\end{rmk}

The inverse of $\delta$ on $\check \Omega = \widetilde{\Ein}^{1,n-1} \setminus \tlambda$ is a conformal embedding $\lambda :  \check \Omega  \to {\widetilde M}$.  Because $n \geq 3$, $\partial \check \Omega$ has  codimension at least $2$.  Then theorem 1.8 of \cite{charles.bords} applies in our context.  
It yields an open subset $\check \Omega^{\prime}$ containing $\check \Omega$ and a conformal diffeomorphism $\lambda^{-1} : {\widetilde M} \to \check \Omega^{\prime}$, which coincides with $\delta$ on $\Omega$.  Two conformal maps which are the same on an open set of a connected pseudo-Riemannian manifold of dimension $\geq 3$ must coincide, so $\lambda^{-1}=\delta$.  Now, $\check \Omega^{\prime}$ contains $\delta(\tdelta) = \tlambda$ and $\widetilde{\Ein}^{1,n-1} \setminus \tlambda$, which yields $\check \Omega^{\prime}= \widetilde{\Ein}^{1,n-1}$, and $\delta : \widetilde{M} \rightarrow \widetilde{\Ein}^{1,n-1}$ is a conformal diffeomorphism.  

Thus $M$ is conformally diffeomorphic to the quotient $\widetilde{\Ein}^{1,n-1}/ \Phi$.  Since $\Phi$ centralizes $\check \h$, it leaves $\tlambda$ invariant by proposition \ref{finite.action}, and the restriction of $\Phi$ to $\tlambda$ factors through a homomorphism to $Z \cong \Z$.  Because $\Phi$ acts freely, this restriction homomorphism is injective.  Since $M$ is compact, $\Phi$ must be infinite, so theorem \ref{flatness.thm} is proved in the Lorentzian case. 

\subsection{Proof of theorem \ref{flatness.thm}: the case $p \geq 2$}
\label{sec.proof.p}
The proof in the Lorentz case must be adapted for $p \geq 2$ because in this case, $\delta$ is {\it a priori} just a covering map from $\tdelta$ to $\tlambda$ and no longer a diffeomorphism.  

Recall that $\tau^s$ fixes $\tlambda$ pointwise.  Let $p_1 \in \tlambda$.  The lightcone $C(p_1)$ has two singular points, $p_1$, and another point, $p_2 \in \tlambda$, and its complement consists of two Minkowski components, $M_1$ and $M_2$.  Also, $\tlambda \setminus \{ p_1, p_2 \} $ has two connected components ${\mathcal I}_1$ and ${\mathcal I}_2$, which can be defined by dynamical properties of $\tau^s$:
\begin{eqnarray*}
\forall z \in M_1, \ \lim_{s \rightarrow \infty} \tau^s.z \in \mathcal{I}_1 & \mbox{ and }  & \lim_{s \rightarrow - \infty} \tau^s.z \in \mathcal{I}_2 \\
\forall z \in M_2, \ \lim_{s \rightarrow \infty} \tau^s.z \in \mathcal{I}_2 & \mbox{ and } & \lim_{s \rightarrow - \infty} \tau^s.z \in \mathcal{I}_1
\end{eqnarray*}

If $F$ is the set of fixed points of $\tau^s$, then $C(p_1) \setminus F$ splits into two connected components, $C_1$ and $C_2$.  Suppose $p_1 = [e_0]$.  Then $C(p_1)$ is the quotient $(e_0^\perp \cap \widehat{\mathcal{N}}^{p+1,q+1})/\R_{>0}^*$ and $p_2 = [-e_0]$.  Recall that 
$$F = (e_0^\perp \cap e_1^\perp \cap \widehat{\mathcal{N}}^{p+1,q+1})/\R_{>0}^*$$
  The components of $C(p_1) \setminus F$ correspond to $\{ \langle x, e_1 \rangle > 0 \}$ and $\{ \langle x, e_1 \rangle < 0 \}$.  As in section \ref{subsection.dynamics.null}, $\tau^s.[x] \rightarrow [e_0]$ as $s \rightarrow \infty$ if $ \langle x , e_1 \rangle > 0$ and $\tau^s.[x] \rightarrow [-e_0]$ if $\langle x, e_1 \rangle < 0$; similarly, $\tau^s.[x] \rightarrow [e_0]$ as $s \rightarrow - \infty$ if $ \langle x , e_1 \rangle < 0$ and $\tau^s.[x] \rightarrow [-e_0]$ as $s \rightarrow - \infty$ if $\langle x, e_1 \rangle > 0$.

The dynamics are the same at any $p_1 \in \tlambda$, because there is a conformal automorphism of $\tein$ sending $[e_0]$ to $p_1$ and preserving $F$: 
\begin{eqnarray*}
\forall z \in C_1, \ \lim_{s \to \infty} \tau^s.z=p_1 & \mbox{ and } & \lim_{s \to - \infty} \tau^s.z=p_2 \\
\forall z \in C_2, \  \lim_{s \to \infty} \tau^s.z=p_2  & \mbox{ and } &  \lim_{s \to - \infty} \tau^s.z=p_1
\end{eqnarray*}

Let $\{\tp_{2i+1} \ : \ i \in J\} = \delta^{-1}(p_1)$ and $\{\tp_{2i} \ : \ i \in J \} = \delta^{-1}(p_2)$.  Order the points $\tp_{2i+1}$ and $\tp_{2i}$ compatibly with an orientation of ${\tilde \Delta}$, and in such a way  that $\tp_{2i}$ is between $\tp_{2i-1}$ and $\tp_{2i+1}$.  If the covering $\delta : {\tilde \Delta} \to {\tilde \Lambda}$ is finite, then $J$ is finite; in this case, order each set of points cyclically.  The segment of ${\tilde \Delta}$ from $\tp_{2i-1}$ to $\tp_{2i}$ will be denoted $I_{2i-1}$, and the segment from $\tp_{2i}$ to $\tp_{2i+1}$ will be $I_{2i}$.  Now the set $\Omega$ of the previous section will be replaced by the two sets
\begin{eqnarray*}
\Omega_1 & = & \{ \tz \in {\widetilde M} \setminus \tdelta \ | \ \lim_{s \to \infty} h^s.\tz \ \mbox{exists and is in } I_1 \} \\
\Omega_2 & = & \{ \tz \in {\widetilde M} \setminus \tdelta \ | \ \lim_{s \to \infty} h^s.\tz \ \mbox{exists and is in } I_2 \}
\end{eqnarray*}

Using the dynamical characterization of $I_1$ and $I_2$ corresponding to that of ${\mathcal I}_1$ and ${\mathcal I}_2$ given above, one can reproduce the proof of proposition \ref{ouverts} to obtain
\begin{propn}
\label{ouverts2}
The sets $\Omega_1$ and $\Omega_2$ are nonempty and open.  Each $\Omega_i$ is mapped diffeomorphically by $\delta$ onto $M_i$, $i = 1,2$.
\end{propn}

\begin{lmm}
\label{boundary.boundary}
$\delta(\partial \Omega_i) \subset C(p_i)$, $i=1,2$.
\end{lmm}

\begin{Pf}
Let $\tz \in \partial \Omega_1$.  By proposition \ref{ouverts2}, $\delta(\tz) \in {\overline M}_1$.  If $\delta(\tz) \in M_1$, the same proposition gives $\tz^{\prime} \in \Omega_1$ such that $\delta(\tz^{\prime})=\delta(\tz)$.  Then if $U^{\prime}$ is a neighborhood of $\tz^{\prime}$ in $\Omega_1$, and if $U$ is a neighborhood of $\tz$ in ${\widetilde M}$, with $U \cap U^{\prime}= \emptyset$, 
$$ \delta(U \cap \Omega_1) \cap \delta(U^{\prime}) \not = \emptyset$$
contradicting the injectivity of $\delta$ on $\Omega_1$.  

The same proof holds if $\tz \in \partial \Omega_2$. \end{Pf}

If $U$ is an open set of a type-$(p,q)$ pseudo-Riemannian manifold $(N,\sigma)$, and if $x \in N$, denote by $C_U(x)$ the set of points in $U$ which can be joined to $x$ by a lightlike geodesic contained in $U$.  

\begin{lmm}
\label{union}
There exists $U$ a neighborhood of $\tp_1$ in ${\widetilde M}$ such that
$$  U \setminus C_U(\tp_1)=(U \cap \Omega_2)  \bigcup (U \cap \Omega_1) $$
\end{lmm}

\begin{Pf}
First choose $U$ a neighborhood of $\tp_1$ that is geodesically convex for some metric in the conformal class, so $U \setminus C_U(\tp_1)$ is a union of exactly two connected components $U_1$ and $U_2$. (Here we use the assumption $p \geq 2$. In the Lorentz case, there would be three connected components.)  We may choose $U$ small enough that $\delta$ maps $U$ diffeomorphically on its image $V$, and $\delta(C_U(\tp_1))=V \cap C(p_1)$. 

First, $U \cap \Omega_1$ and $U \cap \Omega_2$ are both nonempty: remark \ref{adherence} is easily adapted to the current context to show that $I_1 \subset \overline{\Omega_1}$ and  $I_2 \subset \overline{\Omega_2}$.  Assume then that $U_1 \cap \Omega_1 \not = \emptyset$.  By lemma \ref{boundary.boundary}, if $U_1 \cap \partial \Omega_1 \not = \emptyset$, then $\delta(U_1 \cap \partial \Omega_1) \subset C(p_1)$.  
Since $\delta$ is injective on $U$, then $U_1 \cap \partial \Omega_1 \subset C_U(\tilde{p}_1)$, a contradiction. Therefore, $U_1 \cap \partial \Omega_1  = \emptyset$, so $U_1 \subset \Omega_1$.  Similarly, $U_2 \subset \Omega_2$.  \end{Pf}

Let $W_U=C_U(\tp_1) \setminus F$, and define $W=\bigcup_{s \in \R} h^s. W_U$.

\begin{lmm}
The set $\Omega=\Omega_1 \cup W \cup \Omega_2 \subset \widetilde{M}$ is open, and is mapped diffeomorphically by $\delta$ to $\widetilde{\Ein}^{p,q} \setminus F$.
\end{lmm}

\begin{Pf}
We first prove that $\Omega$ is open.  By lemma \ref{union}, and the fact that $\Omega_1$ and $\Omega_2$ are open, the set $\Omega_1 \cup W_U \cup \Omega_2$ is open.  Now, if $\tz \in W$, there exists $s_0 \in \R$ such that $h^{s_0}.\tz \in W_U$.  Then there is a neighborhood $U^{\prime}$ of $h^{s_0}.\tz$ contained in $\Omega_1 \cup W_U \cup \Omega_2 \subset \Omega$.  Then $h^{-s_0}.U^{\prime}$ is a neighborhood of $\tz$ contained in $\Omega$.

We now show that $\delta$ is injective on $\Omega$.  By lemma \ref{ouverts2}, $\delta$ is injective on $\Omega_1$ and $\Omega_2$, and because $\delta(\Omega_1)=M_1$ is disjoint from $\delta(\Omega_2)=M_2$, the map $\delta$ is actually injective on $\Omega_1 \cup \Omega_2$.  Because $\delta(W) \subset C(p_1)$ is disjoint from $\delta(\Omega_1 \cup \Omega_2)$, it suffices to prove that $\delta$ is injective on $W$.  Assume $\tz,\tz^{\prime} \in W$ with $\delta(\tz)=\delta(\tz^{\prime})$, and suppose this point is in $C_1$, so 
$$\lim_{s \to \infty}\tau^s.\delta(\tz)=\lim_{s \to \infty}\tau^s.\delta(\tz^{\prime})=p_1$$

Since $\tz \in W$, either $\lim_{s \to \infty} h^s.\tz=\tp_1$ or $\lim_{s \to - \infty}h^s.\tz=\tp_1$.  But if $\lim_{s \to - \infty}h^s.\tz=\tp_1$, then $\lim_{s \to -\infty}\tau^s.\delta(\tz)=p_1$, contradicting $\delta(\tz) \in C_1$.  Therefore, $\lim_{s \rightarrow \infty} h^s.\tilde{z} = \tilde{p_1}$, and for the same reasons, $\lim_{s \to \infty} h^s.\tz^{\prime}=\tp_1$.  Then there exists $s_0>0$ such that for all $s \geq s_0$, both $h^s.\tz$ and $h^s.\tz^{\prime}$ are in $U$.  Since $\delta(h^s.\tz)=\tau^s.\delta(\tz)=\tau^s.\delta(\tz^{\prime})=\delta(h^s.\tz^{\prime})$, and since $\delta$ is injective on $U$, we get $h^s.\tz=h^s.\tz^{\prime}$ and finally $\tz=\tz^{\prime}$.  The proof is similar if $\delta(\tz) = \delta(\tz')$ is in $C_2$.

It remains to understand the  set $\delta(\Omega)$.  From proposition \ref{ouverts2}, $M_1 \cup M_2 \subset \delta(\Omega)$, and it is also clear that $\delta(\Omega) \subset \tein \setminus F$.  If $z \in C_1$, then there exists $s>0$ such that $\tau^s.z \in V$.  Hence, there is $\tz \in U$ such that $\delta(\tz)=\tau^s.z$, and finally $\delta(h^{-s}.\tz)=z$. Since $\tz \in U$, then $h^{-s}.\tz \in \Omega$, which proves $z \in \delta(\Omega)$. In the same way, we show   that if $z \in C_2$, then $z \in \delta(\Omega)$.  Finally $\delta(\Omega)=\tein \setminus F = M_1 \cup M_2 \cup C_1 \cup C_2$.  
\end{Pf}

The conclusion is essentially the same as in the Lorentzian case.  Let $\check \Omega$ be the complement of $F$ in $\widetilde{\Ein}^{p,q}$.  Then inverting $\delta$ on $\check \Omega$ gives a conformal embedding $\lambda : \check \Omega \to {\widetilde M}$.  Recall from section \ref{subsection.dynamics.null} that $F=\partial \check \Omega$ has codimension $2$.  Then theorem 1.8 of \cite{charles.bords} gives an open subset $\check \Omega^{\prime}$ containing $\check \Omega$ and a conformal diffeomorphism $\lambda^{-1} : {\widetilde M} \to \check \Omega^{\prime}$, which coincides with $\delta$ on $\Omega$.  

As above, $\lambda^{-1}=\delta$.  Now, $\tlambda \subset \check \Omega^{\prime}$, and since $\delta : {\widetilde M} \to \check \Omega^{\prime}$ is a diffeomorphism, the action of $\Phi$ on $\check \Omega^{\prime}$ is free and proper.  In particular, the map associating to an element of $\Phi$ its restriction to $\tlambda$ is injective.  By proposition \ref{finite.action}, the group $\Phi$ is trivial or isomorphic to $\Z_2$.  Because $M = \widetilde{M} / \Gamma$ is compact, $\Phi$ acts cocompactly on $\check \Omega^{\prime}$, so $\check \Omega'= \widetilde{\Ein}^{p,q}$.  Therefore $M$ is conformally diffeomorphic to $\widetilde{\Ein}^{p,q} / \Phi$, as was to be shown.

Charles Frances \\
D\'epartement de Math\'ematiques \\
Universit\'e Paris-Sud XI \\
Charles.Frances@math.u-psud.fr

\smallskip

Karin Melnick \\
(partially supported by NSF fellowship DMS-855735) \\
Department of Mathematics \\
University of Maryland, College Park \\
karin@math.umd.edu

\end{document}